\documentclass[11pt,a4paper]{amsart}
\title[Tropical combinatorics of max-linear Bayesian networks]
      {Tropical combinatorics of max-linear Bayesian networks}
\author{Carlos Am\'endola}
\address{Technical University Berlin, Germany}
\email{amendola@math.tu-berlin.de}
\author{Kamillo Ferry}
\address{Technical University Berlin, Germany}
\email[Corresponding author]{ferry@math.tu-berlin.de}

\date{}


\usepackage{macros}

\usepackage{setspace}
\usepackage[headings]{fullpage}

\usepackage{standalone}

\usetikzlibrary{patterns,patterns.meta}

\usepackage{dsfont}

\usepackage{pgfplots}
\pgfplotsset{compat=1.18}
\usepackage{eqnarray}

\usepackage{faktor}
\usepackage{ytableau}

\usepackage[section]{placeins}

\makeatletter
\AtBeginDocument{%
  \expandafter\renewcommand\expandafter\subsection\expandafter{%
    \expandafter\@fb@secFB\subsection
  }%
}
\makeatother

\usepackage{hyperref}
\hypersetup{
    bookmarksopen,
    bookmarksnumbered,
    pdftitle={Tropical combinatorics of max-linear Bayesian networks},
    pdfauthor={Kamillo Ferry},
    pdfkeywords={project,report},
    pdfstartview={FitH},
    pdfencoding=auto,
    breaklinks=true
}

\usepackage[nameinlink]{cleveref}

\usepackage[normalem]{ulem}

\newtheorem{question}{Question}

\newtheorem*{theorem**}{Theorem\theoremnum}
\newenvironment{theorem*}[1][]{%
  \edef\theoremnum{\if\relax\detokenize{#1}\relax\else~#1\fi}
  \begin{theorem**}
}{%
  \end{theorem**}
}
\newtheorem*{corollary**}{Corollary\theoremnum}
\newenvironment{corollary*}[1][]{%
  \edef\theoremnum{\if\relax\detokenize{#1}\relax\else~#1\fi}
  \begin{corollary**}
}{%
  \end{corollary**}
}

\theoremstyle{definition}
\newtheorem{definition}[theorem]{Definition}
\newtheorem{XxmpX}[theorem]{Example} 
\newenvironment{example}    
  {%
   \pushQED{\qed}\begin{XxmpX}}
  {\popQED\end{XxmpX}}

\theoremstyle{remark}
\newtheorem{remark}[theorem]{Remark}

\newcommand{\Q}{\ensuremath{\mathrm{Q}}}

\newcommand{\trred}{\ensuremath{\mathrm{tr}}}

\newcommand{\localpolytropemoduli}[1][\dag]{\ensuremath{\mathcal{M}^\star_{#1}}}
\newcommand{\closedlocalmoduli}[1][\dag]{\ensuremath{\overline{\localpolytropemoduli[#1]}}}
\newcommand{\polytropemoduli}[1]{\ensuremath{\mathcal{M}^\star_{\nabla,#1}}}

\newcommand{\ray}[1]{\ensuremath{\mathrm{pos}(#1)}}

\newcommand{\tc}{\ensuremath{\mathrm{TC}}}

\newcommand{\Newton}{\ensuremath{\mathcal{N}}}

\newcommand{\defpoly}[1][V]{\ensuremath{\mathcal{Q}_{#1}}}

\newcommand{\Kleene}[1][\dag]{\ensuremath{\mathcal{K}(#1)}}
\newcommand{\sfan}{\ensuremath{\Sigma}}
\newcommand{\subdivision}{\ensuremath{\mathcal{S}}}

\newcommand{\recessioncone}{\ensuremath{\mathrm{rec}}}

\begin{document}

\onehalfspace
\begin{abstract}
    A polytrope is a tropical polyhedron that is also classically convex. We study the tropical combinatorial types
    of polytropes associated to weighted directed acyclic graphs (DAGs). 
    This family of polytropes arises in algebraic statistics when describing
    the model class of max-linear Bayesian networks. We show how the edge weights of a network directly relate to the facet structure of the corresponding polytrope. We also give a classification of polytropes from weighted DAGs at different levels of equivalence. These results give insight on the statistical problem of identifiability for a max-linear Bayesian network.
\end{abstract}

\maketitle

\section{Introduction}
A Bayesian network
is a type of directed graphical model (in the sense of \citeauthorandcite{lauritzen1996graphical}), modeling
causal relations between random variables through a directed acyclic graph (DAG). Max-linear Bayesian networks (MLBNs) were introduced by \citeauthorandcite{Gissibl:2018}, to describe causal relations between variables that take large observed values. A common example is modeling a river network where the extreme values represent possible flooding events \cite{asadi2015extremes,TBK:2024a}.
The starting point of our work is that the structural equations that govern these relations are formulated in the language of tropical geometry \cite{Joswig:ETC,Maclagan-Sturmfels}
and can be described by weighted DAGs \cite{amendola2021markov,amendola2022conditional,tran2022tropical}.

Specifically, let $[n]=\{ 1,2,\dots,n\}$ and consider a weighted DAG \(\dag = ([n],E)\) with edge weights \(w\in\RR_{\geq 0}^E\).
Choosing a labeling of \(\dag\) gives rise to a weighted adjacency matrix
\(C\in\mathbb{R}_{\geq 0}^{n\times n}\), such that \(c_{ij}=w(j\to i)>0\) if and only if \(j\to i\) is in the edge set $E$ of \(\dag\).
Denote $a \vee b = \max \{a, b\}$. Then $G$ describes a \emph{max-linear Bayesian network} \({X = (X_1,\dots,X_n)}\) 
where the random variables~\(X_i\) satisfy the recursive equations 
\begin{equation}\label{eq:recursive-structural-eqn}
    X_i = \bigvee_{j=1}^n c_{ij}X_j \vee Z_i	
\end{equation} 
for \(c_{ij} \geq 0\), \(1\leq i,j\leq n\) and independent random variables \(Z_i\).

Solving the structural equations \eqref{eq:recursive-structural-eqn} purely in terms of the error terms $Z_i$ gives equations 
\begin{equation}\label{eq:max-times-kleene}
    X_i = \bigvee_{j=1}^n c_{ij}^\circledast Z_j
\end{equation} where $c_{ij}^\circledast$ denotes the $(i,j)$ entry of the \emph{max-times Kleene star}
\(C^\circledast := I_n \vee C \vee C^{\circledast 2} \vee \dots \vee C^{\circledast n}\).

Interpreting \eqref{eq:max-times-kleene} as a matrix-vector product and 
applying a coordinate-wise negative logarithmic transformation allows us 
to make use of tropical geometry over the \emph{min-plus semiring} 
\(\mathbb{T} = (\mathbb{R}\cup\{\infty\},\oplus,\odot,\infty,0)\), which is the standard setting in \cite{Joswig:ETC,Maclagan-Sturmfels}.

Similarly to a classical polytope, a tropical polyhedron is a set $P$ which is the 
tropical convex hull of finitely many points. But the tropical convex hull of a point set \(V\)
amounts to taking all possible tropical linear combinations of points in \(V\).
For the tropical polytopes we are interested in, the corresponding matrix is a Kleene star, 
which means that $P$ is not only tropically convex, but also classically convex. 
This kind of polyhedron is called a \emph{polytrope}.
\begin{example}\label{ex:introduction}
    Let $G=\kappa_n$ be the complete directed acyclic graph on \(n\) nodes, 
    that is, \(\kappa_n\) has an edge \(j\to i\) whenever \(j < i\).
    For \(n=3\), \(\kappa_3\) is the \emph{directed triangle} \[
        \begin{tikzcd}[row sep=.5em]
            1 \arrow[rd,"a"] \arrow[dd,"c"] & \\
            & 2 \arrow[ld,"b"] \\
            3
        \end{tikzcd}\quad\text{with lower-triangular weight matrix}\quad C = \begin{pmatrix}
            0 & \infty & \infty \\ 
            a & 0 & \infty \\
            c & b & 0
        \end{pmatrix}.
    \]
    The recursive equations \(X = C\odot X \oplus Z\) associated to this network are \begin{align*}
        X_1 &= Z_1, \\
        X_2 &= aX_1 \oplus Z_2, \\
        X_3 &= cX_1 \oplus bX_2 \oplus Z_3.
    \end{align*} 
    The explicit solution to this recursive system of equations is given by 
    \[
        \begin{pmatrix}X_1 \\ X_2 \\ X_3 \end{pmatrix} =
        \begin{pmatrix}
            Z_1 \\ aZ_1 \oplus Z_2 \\ (c \oplus ab)Z_1 \oplus bZ_2 \oplus Z_3
        \end{pmatrix} = Z_1 \odot \begin{pmatrix}0 \\ a \\ c \oplus ab \end{pmatrix} \oplus Z_2 \odot \begin{pmatrix}\infty \\ 0 \\ b \end{pmatrix} \oplus Z_3 \odot \begin{pmatrix} \infty \\ \infty \\ 0 \end{pmatrix}, 
    \]
    which shows that  \(X\) lies in the tropical polytope spanned by the columns of the (min-plus) Kleene star matrix 
    \[
        C^\star = \begin{pmatrix}
            0 & \infty & \infty \\
            a & 0 & \infty \\
            c\oplus ab & b & 0 
        \end{pmatrix}. \qedhere
    \]  
\end{example}
Initially, much work on polytropes has been done in the bounded case.
\citeauthorandcite{Develin.Sturmfels:2004} defined tropical polytopes and identified polytropes as their basic building blocks, without using the term `polytrope'  yet.
They have been studied under the name \emph{alcoved polytopes} by \citeauthorandcite{Lam:2007} as
polytopes arising from deformations of affine Coxeter arrangements of type \(A_n\). 
\citeauthorandcite{Joswig.Kulas:2010} connected these two points of view and coined the name \emph{polytrope}.

However, the polytropes from max-linear Bayesian networks are necessarily unbounded,
because the corresponding matrix of the weighted DAG is (tropically) lower-triangular. 
More general work covering the unbounded case includes \citeauthorandcite{Fink.Rincon:2015} 
who studied Stiefel tropical linear spaces with arbitrary support.
While polytropes were not their main focus, tropical polyhedra generated by \(n\) vertices and 
\(n\)-dimensional Stiefel tropical linear spaces in dimension \(d\) are related to each other 
as both can be specified in terms of the column space of a matrix, 
and the combinatorics of both objects can be described in terms
of regular subdivisions of subpolytopes of products of simplices \(\Delta_{n-1}\times\Delta_{d-1}\). 
This generalizes the previous results to unbounded polytropes in tropical affine space.

\citeauthor{Joswig.Loho:2016}~\cite{Joswig.Loho:2016} took a different approach and 
characterized the same regular subdivisions as projections from a generalization 
of weighted digraph polyhedra called the \emph{envelope}.
One of our main results concerns the tropical polytope generated by a matrix \(C\), and the
tropical polytope generated by the Kleene star \(C^\star\). 
We are able to characterize how the combinatorial data of the former already determines the combinatorial 
data of the latter. The necessary definitions for this are laid out in \Cref{sec:tropical-preliminaries}.

\begin{theorem*}[\ref{thm:polytrope-cell}]
    Let \(C\in\TT^{n\times n}\) be supported on a DAG \(\dag\). 
    The polytrope \(\tconv(C^\star)\subseteq\tconv(C)\) coincides with the closed cell \(\sector_L(C)\) 
    such that \(L\) is the perfect matching \((1,2,\dots,n)\) in the support graph 
    \(\supportgraph[C]\).
\end{theorem*}

Regarding enumeration of alcoved polytopes, \citeauthor{Joswig.Kulas:2010} classified the combinatorial types 
of alcoved polytopes in dimensions 2 and 3. Then, \citeauthorandcite{Tran:2016} described the combinatorial types 
of alcoved polytopes in terms of cones of a certain Gröbner fan in \(\RR^{n^2-n}\).

On the other hand, \citeauthor{Joswig.Schröter:2022}~\cite{Joswig.Schröter:2022,Joswig.Schröter:2019} 
studied polytropes in the context of \emph{shortest-path problems}. In particular, they studied
the combinatorics of polytropes with polynomial parametrization and gave effective methods to 
compute moduli spaces for polytropes.

In a previous version \cite{Joswig.Schröter:2019}, \citeauthor{Joswig.Schröter:2019}
characterized combinatorial types of polytropes in terms of central subdivisions
of the fundamental polytope \(\fundamentalpolytope\) defined in \Cref{sec:polytropes}. 
Their version of the result focused mainly on the case
where \(\dag\) is the complete directed graph, or in other words, where \(C\) 
has no entry equal to \(\infty\). We adapt their proof to allow for infinite entries, 
giving the following result.
\begin{theorem*}[\ref{thm:central-triangulations}]
    The tropical combinatorial types of full-dimensional polytropes in \(\TTAA^{n-1}\) are in
    bijection with the regular central subdivisions of the fundamental polytope \(\fundamentalpolytope\), where \(G\)
    ranges over transitive directed graphs \(\dag\) on \(n\) nodes.
\end{theorem*}

In \Cref{sec:coarse-classification} we introduce machinery in the form of the 
\emph{weighted transitive reduction} of a weighted DAG, to characterize a minimal facet description for polytropes 
associated to DAGs.
\begin{theorem*}[\ref{thm:transitive-reduction-minimal-facets}]
   The weighted transitive reduction \(\dag^\flat_w\) is the unique minimal weighted DAG 
    such that the weighted digraph polyhedra \(\Q(C)\) and \(\Q(C^\flat)\) coincide.
\end{theorem*}
\Cref{thm:transitive-reduction-minimal-facets} has consequences for statistical inference of the parameters of an MLBN.
Namely, we can only uniquely recover the parameters of an MLBN for edges in the weighted transitive reduction, see \Cref{rem:stats-remark}.

Our main goal is to deeply understand polytropes coming from max-linear Bayesian networks.
This is done through the existing theory of tropical polyhedra and tropical hyperplane arrangements.
From this we obtain \emph{open strata} \(\localpolytropemoduli\) for a moduli space of polytropes 
coming from a fixed DAG \(\dag\) as a polyhedral fan. There is similar work by \citeauthorandcite{OPS:2019}
on valuated matroids with a prescribed basis and the \emph{local Dressians} parametrizing those.

Using the results above, we give a classification of polytropes coming from max-linear Bayesian networks on \(n\) nodes 
in terms of their tropical combinatorial types. We also give a procedure to enumerate generic tropical combinatorial types 
to summarize this classification. Another consequence is a description of the moduli space of such polytropes
in terms of unions of polyhedral fans. There are two natural ways of grouping the corresponding cones, which we also describe.

\subsection*{Outline}
The paper is structured as follows. In \Cref{sec:tropical-preliminaries}, we give the necessary background on
tropical polytopes and their dual subdivisions. 
In \Cref{sec:polytropes}, we recall the specifics related to polytropes and give characterizations of
tropical combinatorial equivalence for polytropes.
In \Cref{sec:coarse-classification}, we turn these characterizations into descriptions of a moduli space 
\(\polytropemoduli{n}\). In \Cref{sec:enumeration} we investigate \(\polytropemoduli{n}\) for small \(n\)
by enumerating generic combinatorial types. We conclude with an outlook including questions guiding future research directions.

\section*{Acknowledgements}
The authors thank Melody Chan, Alex Fink, Anne Frühbis-Krüger, Cécile Gachet, Sami Halaseh, 
Michael Joswig, Francesco Nowell, Yue Ren, Benjamin Schröter and Lena Weis for many helpful discussions and suggestions along the way.
We would also like to thank the two anonymous reviewers for their many helpful comments.

The authors are funded by the Deutsche Forschungsgemeinschaft (DFG, German Research Foundation) under Germany´s Excellence Strategy 
– The Berlin Mathematics Research Center MATH+ (EXC-2046/1, project ID: 390685689), Project AA3-16.

\section{Tropical preliminaries}\label{sec:tropical-preliminaries}
We work mostly over the \emph{min-plus tropical semiring} 
\(\TT = \TT_{\min} = (\mathbb{R}\cup\{\infty\},\oplus,\odot)\)
where \(\oplus = \min\) and \(\odot = +\). Where necessary, we distinguish 
between min-plus and max-plus semiring by decorating with \(\min\) resp.\ \(\max\) 
and using \(\boxplus\) for max-tropical addition.

We denote by \(\TTAA^{d-1}\) the \emph{tropical affine space} which is the 
set of points \(x\in\RR^d\) identified via 
\[(x_1,\dots,x_n)\sim\lambda\odot(x_1,\dots,x_d)=(x_1+\lambda,\dots,x_d+\lambda)\] for all \(\lambda\in\RR\).
This is also called the \emph{tropical projective torus} in the literature. We can also
interpret \(\TTAA^{d-1}\) as the quotient \(\mathbb{R}^d/\mathbb{R}\mathbf{1}\)
of Euclidean space by the subspace spanned by the all-ones vector \(\mathbf{1} = (1,\dots,1)\).
Then, the \emph{tropical projective space} \(\TTPP^{d-1}\) is defined as \[
    \TTPP^{d-1} \coloneq \faktor{\left(\TT^{n}\setminus\{(\infty,\dots,\infty)\}\right)}{\mathbb{R}\mathbf{1}}.
\] The space \(\TTPP^{d-1}\) has a stratification by Euclidean spaces in the following way. 
For an arbitrary subset \(Z\subsetneq[d]\) the \emph{stratum} of \(\TTPP^{d-1}\) given by \(Z\)
is defined as\[
    \TTPP^{d-1}(Z) := \{\,
        x\in\TTPP^{d-1} \mid x_i = \infty \text{ if and only if } i\in Z
    \,\}.
\] We can recover \(\TTAA^{d-1}\) as the stratum \(\TTPP^{d-1}(\varnothing)\) and 
call this the \emph{interior stratum} of \(\TTPP^{d-1}\).
On the other hand if \(Z\neq\varnothing\), then \(\TTPP^{d-1}(Z) \cong \TTAA^{d-1-\#Z}\)
is a \emph{boundary stratum} of \(\TTPP^{d-1}\).

Since our graphs are finite, we write \(\RR^E\) for the set of possible weights on a specific digraph 
\(\dag = (V,E)\) and consider a stratification of the set \(\TT^E\) of weights on \(\dag\) possibly allowing 
\(\infty\). We define the strata indexed by subsets \(E'\subseteq E\)
of the edge sets, or equivalently, by subgraphs \(\subdag\subseteq\dag\) as \[
    \TT^{E}(\subdag) := \{\,
        x\in\TT^E \mid x_{j\to i} < \infty \text{ if and only if } j\to i\in\subdag
    \,\}.
\]

We call a set \(V\) of \(n\) points in \(\TTPP^{d-1}\)  a \emph{tropical point configuration}. 
This can also be interpreted as a matrix \(V\in\TT^{d\times n}\). 

\subsection{Tropical polyhedra in affine space}

We now define the basic objects of consideration. Following \citeauthorandcite{But:2010}, we say that a tropical matrix 
\(V\in\TT^{d\times n}\) is \emph{doubly \(\RR\)-astic}
if every tropical row and column sum is finite. 

\begin{definition}
The \emph{tropical polyhedron} defined by a doubly \(\RR\)-astic matrix \(V\in\TT^{d\times n}\) in \(\TTAA^{d-1}\) 
is the image under left multiplication by $V$, that is, the (min-plus)-linear span of the columns of \(V\):
\[
    \tconv(V) \coloneq \left\{\, 
    V\odot\lambda
    \mathop{\,\bigg\vert\,}
    \lambda\in\TTAA^{n-1}
    \,\right\}\subseteq\TTAA^{d-1}.
\] If \(\tconv(V)\) is bounded, it is called a \emph{tropical polytope}.
\end{definition}
The condition that \(V\) is doubly \(\RR\)-astic ensures that \(\tconv(V)\) is a subset of \(\TTAA^{d-1}\)
instead of living in a boundary stratum of \(\TTPP^{d-1}\).
A set $A$ is called \emph{tropically convex} if it contains for any two points 
\(p,q\in A\) the \emph{tropical line segment} \(\tconv(p,q)\).

Every min-tropical polyhedron can be expressed in terms of max-tropical hyperplanes \cite[Section 6]{Fink.Rincon:2015}. 

\begin{definition}
For any point \(v=(v_1,\dots,v_d)\in\TTPP^{d-1}\), the \emph{tropical hyperplane} \(\mathcal{T}(v)\) defined by \(v\) is the set 
of points \(x\in\TTAA^{d-1}\) such that the maximum in the (max-tropical) linear form \[
    \alpha_v(x) \coloneq -v_1 \odot x_1 \boxplus -v_2\odot x_2 \boxplus \dots \boxplus -v_d \odot x_d.
\] is attained at least twice.
\end{definition}
The choice of negative signs in \(\alpha_v\) makes \(v\) the \emph{apex} of \(\mathcal{T}(v)\).

\begin{figure}[t]
    \centering
    \begin{subfigure}{.32\linewidth}
        \begin{tikzpicture}[
            scale=2.5,
            roundnode/.style={circle, draw=black!85, fill=black!85, thin, scale=0.5},
            pseudonode/.style={circle, draw=black!50, fill=black!50, thin, scale=0.5}
        ]
            \node[roundnode, label={[label distance=1pt]-30:$v$}] (apex) at (0,0) {};
            
            \draw[line width=1pt] (apex) --++ (-1,0); 
            \draw[line width=1pt] (apex) --++ (0,-1); 
            \draw[line width=1pt] (apex) --++ (.70,.70); 

            \node[below left=of apex] {$\sector_1(v)$};
            \node[right=of apex] {$\sector_2(v)$};
            \node[above=of apex] {$\sector_3(v)$};
        \end{tikzpicture}
    \end{subfigure}%
    \begin{subfigure}{.32\linewidth}
        \begin{tikzpicture}[
            scale=2.5,
            roundnode/.style={circle, draw=black!85, fill=black!85, thin, scale=0.5},
            pseudonode/.style={circle, draw=black!50, fill=black!50, thin, scale=0.5}
        ]
            \node[roundnode, label={[label distance=1pt]-30:$v = (\infty,\ast,\ast)$}] (apex) at (0,0) {};

            \node at (.8,.4) {$\sector_2(v)$};
            \node at (.4,.8) {$\sector_3(v)$};

            \draw[line width=1pt] (apex) edge[dashed]++ (0.2,0.2) ++ (0.2,0.2) --++ (1.2,1.2); 
        \end{tikzpicture}
    \end{subfigure}%
    \begin{subfigure}{.32\linewidth}
        \begin{tikzpicture}[
            scale=2.2,
            roundnode/.style={circle, draw=black!85, fill=black!85, thin, scale=0.5},
            pseudonode/.style={circle, draw=black!50, fill=black!50, thin, scale=0.5}
        ]
            \node[roundnode, label={[label distance=1pt]-90:$v = (\ast,\infty,\ast)$}] (apex) at (0,0) {};

            \node at (0,-1) {};

            \node at (-1,-.4) {$\sector_1(v)$};
            \node at (-1,.4) {$\sector_3(v)$};
            
            \draw[line width=1pt] (apex) edge[dashed]++ (-0.3,0) ++ (-0.3,0) --++ (-1.7,0); 
        \end{tikzpicture}
    \end{subfigure}
    \caption{Max-tropical hyperplanes \(\mathcal{T}(v)\) in \(\TTAA^2\)
    and their cells. In the last two examples, the apices are at infinity.}
    \label{fig:basic_trop_halfspaces}
\end{figure}

For any \(i\in[d]\) with \(v_i \neq \infty\), the \emph{\(i\)-th open (max-)cell} of \(\mathcal{T}(v)\) is the set
\[
    \sector_i^\circ(v) \coloneq \left\{\, 
        x\in\TTAA^{d-1} 
        \mid 
        -v_i \odot x_i  > -v_\ell \odot x_\ell \,\text{ for all } 1\leq \ell\leq d, i\neq\ell 
    \,\right\},
\] that is, the set of points \(x\in\TTAA^{d-1}\) for which \(-v_i\odot x_i\) realizes the maximum in \(\alpha_v(x)\) uniquely.
The \emph{\(i\)-th closed (max-)cell} \(\sector_i(v)\) of \(\mathcal{T}(v)\) is 
the set of points \(x\in\TTAA^{d-1}\) for which \(-v_i\odot x_i\) realizes the maximum in \(\alpha_v(x)\).
Any hyperplane \(\mathcal{T}(v)\) divides \(\TTAA^{d-1}\) into up to \(d\) sectors each with apex at \(v\), 
one for each monomial of \(\alpha_v\), as seen in \Cref{fig:basic_trop_halfspaces}.

For any point configuration \(V = \{\, v^{(1)}, \dots, v^{(n)} \,\}\subset\TTPP^{d-1}\), the union of the hyperplanes 
\(\mathcal{T}_j \coloneq\mathcal{T}(v^{(j)})\) gives a \emph{tropical hyperplane arrangement} \(\mathcal{T}(V)\)
in \(\TTAA^{d-1}\). Its \emph{defining polynomial} \(\defpoly\) is obtained 
by multiplying the linear forms \(\alpha_{v^{(j)}}(x)\) for the tropical hyperplanes \(\mathcal{T}_j\) 
where \(1\leq j \leq n\). That is, \[
    \defpoly \defeq \bigodot_{j=1}^n \alpha_{v^{(j)}}(x).
\]
\begin{example}\label{ex:tropical-hyperplane-arrangement}
    The \emph{directed triangle} \(\kappa_3\) from \Cref{ex:introduction}
    describes a point configuration $V$ given by the columns of the weighted adjacency matrix \[
        C = \begin{pmatrix}
            0 & \infty & \infty \\ 
            a & 0 & \infty \\
            c & b & 0
        \end{pmatrix}.
    \]  
    This point configuration gives rise to the tropical hyperplane arrangement \(\mathcal{T}(V)\) 
    with defining tropical polynomial \begin{align*}
        \defpoly
        &= \alpha_{v^{(1)}} \odot \alpha_{v^{(2)}} \odot \alpha_{v^{(3)}} \\ 
        &= (x_1 \boxplus -a \odot x_2 \boxplus -c\odot x_3) \odot (x_2 \boxplus -b\odot x_3) \odot x_3 \\
        &= x_1x_2x_3 \boxplus -a \odot x_2^2x_3 \boxplus -(c \oplus ab) \odot x_2x_3^2 \boxplus -b \odot x_1x_3^2 \boxplus -(bc)\odot x_3^3.
    \end{align*} 
    The tropical hyperplane arrangement for \(c < a+b\) is sketched in \Cref{fig:hyperplane-arrangement-sectors}
    with \(\tconv(V)\) shown as the shaded region.
    \begin{figure}[t]
        \centering
        \begin{tikzpicture}[
    scale=2.5,
    roundnode/.style={circle, draw=black!85, fill=black!85, thin, scale=0.5},
    pseudonode/.style={circle, draw=black!85, fill=white, scale=0.5}
]
    \node at (.5,.5) {};
    \node at (0,-1) {};

    \coordinate (v1) at (0,0);
    \coordinate (v2) at (-1,0);

    \coordinate (v3) at (-2,-1);
    \coordinate (v4) at (0,-1);

    \coordinate (u1) at (.5,0);
    \coordinate (u2) at (-2,-1);

    \coordinate (u3) at (0,.5);

    \draw[black!15, fill=black!15] (v1) -- (v2) -- (v3) --++ (2,0) -- cycle;

    \draw[black!85!black, line width=2pt] (v2) -- (v3); 
    \draw[black!85!black, line width=2pt, loosely dotted] (v2) ++ (.5,.5) -- (v2); 
    \draw[black!85!black, line width=2pt] (v2) -- (v1);
    \draw[black!85!black, line width=2pt, loosely dotted] (v2) ++ (-1,0) -- (v2);
    \draw[black!85!black, line width=2pt] (v1) -- node[near end,inner sep=0] (s0) {}  (v4);
    \draw[black!85!black, line width=2pt, loosely dotted] (v1) ++ (.5,.5) edge node[inner sep=0] (s1) {} (v1);

    \node[roundnode] at (v1) {};
    \node[pseudonode] at (v2) {};

    \node[right of=v1,xshift=1cm,inner sep=1pt] (t1) {$\scriptstyle (1,12,13)$};
    \draw [-latex,thick] (t1) to[bend left=15] (v1);

    \node[left of=v2,xshift=-.5cm,yshift=.5cm,inner sep=1pt] (t2) {$\scriptstyle (1,2,123)$};
    \draw [-latex,thick] (t2) to[bend left=15] (v2);

    \node[right of=s1,xshift=1cm,inner sep=1pt] (t3) {$\scriptstyle (\varnothing,12,13)$};
    \draw [-latex,thick] (t3) to[bend left=15] (s1);
    
    \node[right of=s0,xshift=1cm,inner sep=1pt] (t0) {$\scriptstyle (1,12,3)$};
    \draw [-latex,thick] (t0) to[bend left=15] (s0);

    \node[black!85!black, below =of v2, anchor=center, xshift=.5cm] {$(1,2,3)$};
    \node[black!85!black, below left=of v2, anchor=center, xshift=-1.5cm] {$(1,\varnothing,23)$};
    \node[black!85!black, below right=of v1, anchor=center] {$(\varnothing,12, 3)$};
    \node[black!85!black, above left=of v1, anchor=south west] {$(\varnothing, 2, 13)$};
    \node[black!85!black, above left=of v2, xshift=-1cm, anchor=south west] {$(\varnothing, \varnothing, 123)$};
\end{tikzpicture}
        \caption{Max-tropical hyperplane arrangements with some affine covectors.
        The shaded region is \(\tconv(V)\) and the solid dot is the tropical vertex~\(v^{(1)}\).}
        \label{fig:hyperplane-arrangement-sectors}
    \end{figure}
\end{example}

The polyhedral subdivision of \(\TTAA^{d-1}\) induced by \(\mathcal{T}(V)\) 
is called the \emph{covector decomposition} \(\CovDec(V)\) induced by \(V\) and it
is the common refinement of the subdivisions for each \(\mathcal{T}_j\) \cite{Fink.Rincon:2015}. 
The \emph{covector cells} of \(\CovDec(V)\) may be described combinatorially in the following way.

Define the \emph{support graph} \(\supportgraph\) of \(V\) as the bipartite graph on \([n]\amalg[d]\) 
with edges \[\left\{\, (j,i) \in [n]\times[d] \mid v_i^{(j)} \neq \infty \,\right\}.\]
A \emph{\((d,n)\)-covector} \(L\) is a subgraph of the support graph \(\supportgraph\).
If clear from the context, we might refer to \(L\) as just a \emph{covector}.
For a covector \(L\), write \(L_i\) for the predecessors of the node \(i\in[d]\) in \(L\).
This also allows for a compact notation of a covector as \(L = (L_1,\dots,L_d)\), compare with 
\cite[Section 6.3.]{Joswig:ETC}. For example, in compact notation, \((1,\varnothing,23)\) refers to
the covector \[
    \begin{tikzcd}[font = \small, row sep=1.75em, column sep=2em]
        1 \arrow[d] & 2 \arrow[dr] & 3 \arrow[d] \\
        1 & 2 & 3
    \end{tikzcd}.
\]
Any polyhedral cell of \(\CovDec(V)\) can be labeled in the following way. 

\begin{definition}
The \emph{affine covector} of a point \(p\in\TTAA^{d-1}\) relative to \(V\) 
is the covector \(\Cov(p,V)\) containing edges \((j,i)\) if and only if 
\(p\) is contained in the sector \(\sector_i(v^{(j)})\) of the vertex \(v^{(j)}\).
\end{definition}

The \emph{open (covector) cell} \(\sector^\circ_L(V)\) is the intersection of the open cells 
\(\sector_i^\circ(v^{(j)})\) for every edge \((j,i)\in L\).
Similarly, the \emph{closed (covector) cell} \(\sector_L(V)\) is the intersection of the closed cells 
\(\sector_i(v^{(j)})\) for all \((j,i)\in L\). 

Both inclusion \(L\subseteq L'\) and union \(L\cup L'\) of \((d,n)\)-covectors are to be understood as 
the corresponding operations on the edge sets.
From the definitions it follows that
\begin{equation}
    \sector_{L\cup L'}(V) = \sector_L(V) \cap \sector_{L'}(V). \label{eq:cap-sectors-cup-covector}
\end{equation} 
In particular, the full-dimensional open cells 
correspond to covectors where every vertex \(j\in [n]\) only has one outgoing edge.

\begin{remark}\label{rem:sectors-monomials}
Each full-dimensional open cell corresponds to a monomial \(x^a = x_1^{a_1}\dots{}x_d^{a_d}\) 
in the following way. Let \(S\) be a full-dimensional open cell of a covector decomposition \(\CovDec\).
A generic point \(x\in S\) has a covector \(\Cov(x,V) = (L_1,\dots,L_d)\) which encodes the fact that \[
    \alpha_j(x) = -v_{ij}\odot x_i
\] whenever \(j\in L_i\). Since \(\defpoly = \odot_{j=1}^n \alpha_j\), this implies that for a generic point \(x\in S\)
that \[
    \defpoly(x) = \bigodot_{j=1}^n (-v_{ij}\odot x_i)
\] which is exactly the monomial of \(\defpoly\) with \(x_i^{a_i}\) for \(a_i = \#L_i\).
In the context of tropical hyperplane arrangements, \citeauthor{Dochtermann.Joswig.ea:2012}
referred to the vector of counts \(a = (\#L_1,\dots,\#L_d)\) as \emph{coarse type} of the point \(x\) 
and made the connection between coarse types and exponent vectors of monomials in \(\defpoly\) \cite{Dochtermann.Joswig.ea:2012}.

Then, the relationship between open sectors and covectors in \eqref{eq:cap-sectors-cup-covector} 
implies that the corresponding coarse type of a lower-dimensional open sector \(S\) does not correspond 
to a monomial of \(\defpoly\). Instead, if \(S = \cap_{i\in I} S_i\) where \(S_i\) are full-dimensional closed sectors,
then the coarse type of \(S\) is the least common multiple of the coarse types of each \(S_i\).
\end{remark}

Note that while the correspondence between full-dimensional open cells and monomials is one-to-one,
there are possibly multiple covectors leading to the same monomial of \(\defpoly\). This will be 
discussed in more detail in the next section.

\begin{example}\label{ex:affine-covectors}
    Continuing from \Cref{ex:tropical-hyperplane-arrangement}, 
    \Cref{fig:hyperplane-arrangement-sectors} contains a non-exhaustive set of examples 
    for affine covectors relative to the point configuration \(V\).
    The support graph \(\supportgraph[V]\) of this point configuration is given by \[
        \begin{tikzcd}[font = \small, row sep=1.75em, column sep=2em]
            1 \arrow[d] \arrow[dr] \arrow[drr] & 2 \arrow[d]\arrow[dr] & 3 \arrow[d] \\
            1 & 2 & 3
        \end{tikzcd}.
    \] With respect to \Cref{rem:sectors-monomials}, there are two covectors giving rise to the monomial \(x_2x_3^2\),
    namely \[
        \begin{tikzcd}[font = \small, row sep=1.75em, column sep=2em]
            1 \arrow[drr,pos=.35,swap,"c"] & 2 \arrow[d,pos=.2,"0"] & 3 \arrow[d,pos=.2,"0"] \\
            1 & 2 & 3
        \end{tikzcd}
        \quad\text{and}\quad
        \begin{tikzcd}[font = \small, row sep=1.75em, column sep=2em]
            1 \arrow[dr,pos=.35,swap,"a"] & 2 \arrow[dr,pos=.35,swap,"b"] & 3 \arrow[d,pos=.35,"0"] \\
            1 & 2 & 3
        \end{tikzcd}.
    \] The edges \(j\to i\) of both covectors are annotated with the values for \(c_{ij}\).
    Multiplying together the terms of each \(\alpha_j\) corresponding in both covectors
    gives \[
        \left(-(c \odot 0 \odot 0) \boxplus -(a\odot b \odot 0)\right) x_2x_3^2
        = -(c\oplus ab) x_2x_3^2.
    \] This suggests that the coefficient of \(x_2x_3^2\) arises as minus the minimum over the 
    sums of weights for covectors \((L_1,L_2,L_3)\) with \(a_i = \#L_i\).
    Similarly, all five maximal-dimensional cells correspond to the five monomials of 
    the defining polynomial \(\defpoly[V]\). \qedhere
\end{example}

It follows from \eqref{eq:cap-sectors-cup-covector} that the poset of affine covectors of 
a point configuration \(V\subset\TTPP^{d-1}\) is dual to the poset of cells in \(\CovDec(V)\). 
This justifies the following definition, the terminology following  \citeauthor{Fink.Rincon:2015}~\cite{Fink.Rincon:2015}.

\begin{definition}
    The \emph{tropical combinatorial type} \(\tc(V)\) of a point configuration \(V\) is defined
    as the poset of affine covectors of $V$.
\end{definition}

A min-tropical polytope \(\tconv(V)\) inherits its tropical combinatorial type from 
the max-tropical hyperplane arrangement \(\mathcal{T}(V)\) in the following way.
\begin{proposition}[{\cite[Thm.\ 6.\,14.]{Joswig:ETC}}]\label{thm:covdec-tconv-cells}
    For \(V\subset\TTPP^{d-1}\) finite, the tropical polyhedron \(\tconv(V)\) in \(\TTAA^{d-1}\) is the union 
    of closed cells with covector \(L\) such that every \(L_i \neq \varnothing\) and~\(\bigcup_{i=1}^d L_i = [n]\).
    In other words, every node in \(L\) has an incoming and outgoing edge.
\end{proposition}

\begin{example}
    The covectors for the arrangement \(\mathcal{T}(V)\) from \Cref{ex:affine-covectors} such that no node in the $(3,3)$-covector
    is isolated are \[
        (1,2,3),\ (1,12,3),\ (1,2,23),\ (1,2,13),\ (1,12,13),\ \text{and}\ (1,2,123).
    \] These are precisely the covectors associated to the closed shaded region in 
    \Cref{fig:hyperplane-arrangement-sectors}.
\end{example}

If the tropical vertex set \(V\subset\TTAA^{d-1}\), that is, every vertex has only finite
coordinates, the above statement would admit the slogan version `\textit{\(\tconv(V)\) 
is the union of the bounded cells in \(\CovDec(V)\)}' which is an important result by 
\citeauthor{Develin.Sturmfels:2004}~\cite[Theorem 15]{Develin.Sturmfels:2004}.

\begin{definition}\label{def:tropically-equivalent}
    Two point configurations \(V = (v_1,\dots,v_n)\) and \(W = (w_1,\dots,w_n)\) of points in \(\TTPP^{d-1}\) 
    are \emph{tropically equivalent} if there is a pair of permutations 
    \((\tau, \sigma) \in \SymGr_n \times \SymGr_d\) acting on the vertices and coordinates
    such that the induced map on covectors is a poset isomorphism between \(\tc(V)\) and \(\tc(W)\).
\end{definition}

Similarly, two tropical polyhedra \(\tconv(V)\) and \(\tconv(W)\) are \emph{tropically equivalent} if the point configurations of \(V\) and \(W\) are tropically equivalent.

\subsection{Regular dual subdivisions}\label{sec:regular-dual-subdivisions}
A polyhedral subdivision of a finite point set \(A\) in \(\RR^d\) is called \emph{regular} if it is
induced by a height function \(h\in\RR^A\). We denote this subdivision by \(\subdivision(A,h)\).
For further details, see \cite[Section A.\,4.]{Joswig:ETC}.

\begin{definition}
Let \(\subdivision\) be a subdivision of a point configuration \(A\subset\RR^{d}\)
and consider all regular subdivisions \(\subdivision(A,h)\) that are refined by \(\subdivision\).

The associated height functions form the \emph{secondary cone} \[
    \sfan(A,\subdivision) = \left\{\, 
        h\in\RR^A \mid \subdivision \text{ refines } \subdivision(A,h) 
    \,\right\}.
\] The collection of all secondary cones forms the \emph{secondary fan} \(\sfan(A)\).
\end{definition}

\begin{example}\label{ex:basic-example-sfan}
    \begin{figure}[b]





\hspace*{\fill}
\begin{tikzpicture}[
    scale=2.5,
    roundnode/.style={circle, draw=black!85, fill=black!85, thin},
    pseudonode/.style={circle, draw=black!50, fill=black!50},
    x = 1em, y = 1em
]
    \node[roundnode, inner sep=1pt, label={[label distance=1em]-90:$0$}] (o) at (0,0) {};
    \node[inner sep=1pt, label={[label distance=1em]-90:$1$}] (i) at (1,0) {};
    \node[roundnode, inner sep=1pt, label={[label distance=1em]-90:$2$}] (ii) at (2,0) {};
    \node[roundnode, inner sep=1pt, label={[label distance=1em]-90:$3$}] (iii) at (3,0) {};

    \node[roundnode, inner sep=1pt, draw=blue!85, fill=blue!85] (ho) at (0,1) {};
    \node[pseudonode, inner sep=1pt] (hi) at (1,2) {};
    \node[roundnode, inner sep=1pt, draw=blue!85, fill=blue!85] (hii) at (2,.5) {};
    \node[roundnode, inner sep=1pt, draw=blue!85, fill=blue!85] (hiii) at (3,2) {};

    \draw[black] (o) -- (iii);
    \draw[blue] (ho) -- (hii) -- (hiii);
    \draw[blue!50] (ho) --++ (0,2);
    \draw[blue!50] (hiii) --++ (0,1);
\end{tikzpicture}
\hfill
\begin{tikzpicture}[
    scale=2.5,
    roundnode/.style={circle, draw=black!85, fill=black!85, thin},
    pseudonode/.style={circle, draw=black!50, fill=black!50},
    x = 1em, y = 1em
]
    \node[roundnode, inner sep=1pt, label={[label distance=1em]-90:$0$}] (o) at (0,0) {};
    \node[roundnode, inner sep=1pt, label={[label distance=1em]-90:$1$}] (i) at (1,0) {};
    \node[inner sep=1pt, label={[label distance=1em]-90:$2$}] (ii) at (2,0) {};
    \node[roundnode, inner sep=1pt, label={[label distance=1em]-90:$3$}] (iii) at (3,0) {};

    \node[roundnode, inner sep=1pt, draw=blue!85, fill=blue!85] (ho) at (0,1) {};
    \node[roundnode, inner sep=1pt, draw=blue!85, fill=blue!85] (hi) at (1,.5) {};
    \node[pseudonode, inner sep=1pt] (hii) at (2,2) {};
    \node[roundnode, inner sep=1pt, draw=blue!85, fill=blue!85] (hiii) at (3,2) {};

    \draw[black] (o) -- (iii);
    \draw[blue] (ho) -- (hi) -- (hiii);
    \draw[blue!50] (ho) --++ (0,2);
    \draw[blue!50] (hiii) --++ (0,1);
\end{tikzpicture}
\hfill
\begin{tikzpicture}[
    scale=2.5,
    roundnode/.style={circle, draw=black!85, fill=black!85, thin},
    pseudonode/.style={circle, draw=black!50, fill=black!50},
    x = 1em, y = 1em
]
    \node[roundnode, inner sep=1pt, label={[label distance=1em]-90:$0$}] (o) at (0,0) {};
    \node[roundnode, inner sep=1pt, label={[label distance=1em]-90:$1$}] (i) at (1,0) {};
    \node[roundnode, inner sep=1pt, label={[label distance=1em]-90:$2$}] (ii) at (2,0) {};
    \node[roundnode, inner sep=1pt, label={[label distance=1em]-90:$3$}] (iii) at (3,0) {};

    \node[roundnode, inner sep=1pt, draw=blue!85, fill=blue!85] (ho) at (0,1) {};
    \node[roundnode, inner sep=1pt, draw=blue!85, fill=blue!85] (hi) at (1,.5) {};
    \node[roundnode, inner sep=1pt, draw=blue!85, fill=blue!85] (hii) at (2,1) {};
    \node[roundnode, inner sep=1pt, draw=blue!85, fill=blue!85] (hiii) at (3,2) {};

    \draw[black] (o) -- (iii);
    \draw[blue] (ho) -- (hi) -- (hii) -- (hiii);
    \draw[blue!50] (ho) --++ (0,2);
    \draw[blue!50] (hiii) --++ (0,1);
\end{tikzpicture}
\hfill
\begin{tikzpicture}[
    scale=2.5,
    roundnode/.style={circle, draw=black!85, fill=black!85, thin},
    pseudonode/.style={circle, draw=black!50, fill=black!50},
    x = 1em, y = 1em
]
    \node[roundnode, inner sep=1pt, label={[label distance=1em]-90:$0$}] (o) at (0,0) {};
    \node[inner sep=1pt, label={[label distance=1em]-90:$1$}] (i) at (1,0) {};
    \node[inner sep=1pt, label={[label distance=1em]-90:$2$}] (ii) at (2,0) {};
    \node[roundnode, inner sep=1pt, label={[label distance=1em]-90:$3$}] (iii) at (3,0) {};

    \node[roundnode, inner sep=1pt, draw=blue!85, fill=blue!85] (ho) at (0,1) {};
    \node[pseudonode, inner sep=1pt] (hi) at (1,1.5) {};
    \node[pseudonode, inner sep=1pt] (hii) at (2,2) {};
    \node[roundnode, inner sep=1pt, draw=blue!85, fill=blue!85] (hiii) at (3,1) {};

    \draw[black] (o) -- (iii);
    \draw[blue] (ho) -- (hiii);
    \draw[blue!50] (ho) --++ (0,2);
    \draw[blue!50] (hiii) --++ (0,2);
\end{tikzpicture}
\hspace*{\fill}
        \caption{The four possible regular subdivisions of the point configuration 
        \(A = \{\, 0,1,2,3\,\}\subset\RR\) in \Cref{ex:basic-example-sfan}.}
        \label{fig:basic-example-sfan}
    \end{figure}
    Consider the point configuration \(A = \{\, 0,1,2,3\,\}\subset\RR\). 
    This has four possible regular subdivisions shown in \Cref{fig:basic-example-sfan}. 
    Its secondary fan \(\sfan(A)\) in \(\RR^4\) is a complete 
    (\ie{} covering the entire space) full-dimensional polyhedral fan. 
    It has four cones corresponding to the four possible regular subdivisions. 
\end{example}

Denote by \(\supp(F)\) the support of a tropical polynomial \(F \in \TT[x_1,\dots,x_d]\), i.e. the set of exponent vectors 
\(s\in\mathbb{Z}^d\) such that the coefficient \(c_s \neq \infty\).
The \emph{Newton polytope} \(\Newton(F)\) of a tropical polynomial \(F\) is the convex hull of \(\supp(F)\) in \(\mathbb{R}^d\).

Note, that if \(V\subset\TT\mathbb{A}^{d-1}\), then \(\supp(\defpoly) = n\cdot \Delta_{d-1}\),
as \(\defpoly\) is a product of linear forms with Newton polytope \(\Delta_{d-1}\) each.
In the general case, \citeauthor{Fink.Rincon:2015}~\cite{Fink.Rincon:2015} gave the following description 
of \(\Newton(\defpoly)\) in terms of \emph{mixed subdivisions}. For details we refer to \cite[Section 9.\,2.\,2.]{DRS:2010}.

\begin{proposition}[{\cite[Proposition 4.\,1]{Fink.Rincon:2015}}]\label{prop:mixed-regular-subdivision-newton}
  The covector decomposition \(\CovDec(V)\) is dual to the mixed regular subdivision of \(\Newton(\defpoly)\)
  with height function given by the coefficients of \(\defpoly\).
  A face of \(\CovDec(V)\) labeled by an affine covector \(L\) is dual to the cell of \(\Newton(\defpoly)\)
  obtained as the Minkowski sum 
  \begin{equation}\label{eq:minkowski-sum-of-simplices}
      \Newton(\defpoly) = \sum_{j=1}^n \Delta_{L_i}
  \end{equation} where \(\Delta_{L_i}\) denotes the simplex with vertices \(e_j\) for \(j\in L_i\).
\end{proposition}

Instead of enumerating every regular subdivision of \(\Newton(\defpoly)\)
and checking whether a given subdivision is actually mixed, we can use the
\emph{Cayley trick} (see \cite{santos2005cayley}) to check regular subdivisions of \(\Delta_{n-1}\times\Delta_{d-1}\) and its subpolytopes instead. 

\begin{definition}
    Any subpolytope of \(\Delta_{n-1}\times\Delta_{d-1}\) is called a \emph{root polytope}. The \emph{root polytope} \(\rootpolytope[V]\)
    associated to point configuration \(V\in\TT^{d\times n}\) is given by \[
        \rootpolytope[V]\coloneq \conv\left\{\,
            (e_i, e_j) \mid v_i^{(j)} < \infty
        \,\right\}.
    \]
\end{definition}

Through the Cayley trick we get a bijection between subdivisions of \(\rootpolytope[V]\) and mixed subdivisions 
of \(\Newton(\defpoly)\). The following result has been termed the \emph{structure theorem of tropical convexity}
by \citeauthorandcite{Joswig.Loho:2016}.

\begin{theorem}[{\cite[Cor.\ 34]{Joswig.Loho:2016}}]\label{thm:structure-theorem}
    The covector decomposition \(\CovDec(V)\) of \(\TTAA^{d-1}\) is dual to the regular subdivision 
    \(\subdivision\) of the root polytope \(\rootpolytope[V]\) with height function given by \(V\).

    Moreover, the covector decomposition of the tropical polytope \(\tconv(V)\) is dual to 
    \(\subdivision\) intersected with the interior of \(\Delta_{n-1}\times\Delta_{d-1}\).
\end{theorem}

\begin{corollary}
    There is a bijection between combinatorial types of tropical polytopes 
    with \(n\) vertices in \(\TTAA^{d-1}\) and regular subdivisions of root polytopes in 
    \({\Delta_{n-1}\times\Delta_{d-1}}\).
\end{corollary}

In what follows, we specialize this result to the case of polytropes 
to later enumerate combinatorial types of polytropes in Section \ref{sec:enumeration}.

\subsection{Compactifications of polyhedral fans}\label{sec:compactifications}
We recall the definitions of polyhedral spaces and compactifications of polyhedra in tropical toric varieties.
The definitions follow \citeauthorandcite{KSW:2023}.

\begin{definition}
    For a rational polyhedral fan \(\Delta\subset N_\RR\) the \emph{tropical toric variety} \(N_\RR(\Delta)\) of \(N_\RR\)
    with respect to \(\Delta\) is \[
        N_\RR(\Delta) \defeq \coprod_{\sigma\in\Delta} N_\RR/\mathrm{span}(\sigma).
    \] where \(N_\RR(\sigma) \coloneqq  N_\RR/\mathrm{span}(\sigma)\) is the \emph{stratum} associated to \(\sigma\).
    We equip \(N_\RR(\Delta)\) with the unique topology such that \begin{itemize}
        \item the inclusions \(N_\RR/\mathrm{span}(\sigma)\hookrightarrow N_\RR(\Delta)\) are continuous for any cone \(\sigma\in\Delta\), and
        \item for any \(x,v\in N_\RR\), the sequence \((x+nv)_{n\in\NN}\) converges in \(N_\RR(\Delta)\) if and only if 
            \(v\) is contained in the support of the fan \(\Delta\).
    \end{itemize}
\end{definition}

For any cone \(\sigma\in\Delta\) we denote by \(\pi_\sigma\) the natural projection maps \(N_\RR\twoheadrightarrow N_\RR(\sigma)\) onto each stratum.

\begin{remark}
    The tropical toric variety we are mostly interested in is given by the cone \(\Delta\) over the 
    \((d-1)\)-simplex \(\Delta_{d-1}\). In this case, \(N_\RR(\Delta) = \TT^d\).
    Similarly, in the case where we take the cone \(\Delta\) inside \(\RR^E\) spanned by all positive 
    coordinate directions, we get \(N_\RR(\Delta) = \TT^E\). Note that in this case, we index the strata of \(\TT^E\)
    by the coordinates that are not equal to \(\infty\), as opposed to the charts of \(N_\RR(\Delta)\) where
    a cone \(\sigma\) denotes the coordinates that are set to infinity.
\end{remark}

\begin{definition}[{\cite[Def 3.1]{OR:2013}}]
  Let \(\mathcal{P}\) be a finite collection of polyhedra in \(N_\RR\), and \(\Delta\) be a pointed fan.
  The fan \(\Delta\) is said to be \emph{compatible} with \(\mathcal{P}\) if for all \(P\in\mathcal{P}\)
  and all cones \(\sigma\in\Delta\), either \(\sigma\subset\relint(P)\) or \(\relint(\sigma)\cap\recessioncone(P) = \varnothing\), where $\recessioncone(P)$ denotes the recession cone of $P$. 
If for all \(P\in\mathcal{P}\), $\recessioncone(P)$ is a union of cones in \(\Delta\),
  we call \(\Delta\) a \emph{compactifying fan} for \(\mathcal{P}\).
\end{definition}

A special case for compactifying fans is the compactification of a single polyhedron \(P\). In this case, 
reasonable choices of compactifying fans are given by faces of the recession cone. This leaves us with the following description
of the compactification.

\begin{lemma}[{\cite[Lemma 3.9]{OR:2013}}]\label{lem:compactifying-fan}
    Let \(\Delta\) be a compactifying fan for the polyhedron \(P\). 
    Then, its compactification \(\overline{P}\) in \(N_\RR(\Delta)\)
    is \[
        \overline{P} \defeq \coprod_{\sigma\in\Delta, \relint(\sigma)\cap\rec(P)\neq\varnothing} \pi_\sigma(P).
    \] 
\end{lemma}

Recall from \citeauthorandcite{JSS:2019} that a polyhedral subspace 
is the
underlying set of a polyhedral complex in \(\TT^n\). The definition of a polyhedral complex in \(\TT^n\)
is analogous to the definition of a polyhedral complex in Euclidean space, where the comprising polyhedra in \(\TT^n\)
are taken as closures of polyhedra in the strata of \(\TT^n\) (with respect to the topology as a tropical toric variety).
Thus we can make the following definition, which we include for completeness.

\begin{definition}
    A \emph{polyhedral space} \(X\) is a paracompact, second-countable Hausdorff space with an atlas of charts 
    \((\varphi_\alpha\colon U_\alpha\to\Omega_\alpha\subset X_\alpha)_{\alpha\in A}\) such that:
    \begin{enumerate}
        \item The \(U_\alpha\) are open subsets of \(X\), the \(\Omega_\alpha\) are open subsets of polyhedral subspaces
            \(X_\alpha \subset \TT^{r_\alpha}\), and the maps \(\varphi_\alpha\colon U_\alpha\to\Omega_\alpha\) 
            are homeomorphisms for all \(\alpha\in A\), and
        \item for all \(\alpha,\beta\in A\), the transition maps \[
            \varphi_\alpha\circ\varphi_\beta^{-1}\colon\varphi_\beta(U_\alpha\cap U_\beta)\to\varphi(U_\alpha\cap U_\beta)
        \] are extended affine linear maps in the sense of \citeauthorandcite{JSS:2019}.
    \end{enumerate}
\end{definition}

\section{Polytropes from DAGs}\label{sec:polytropes}
We now shift our focus to point configurations \(C\in\TT^{n\times n}\) arising from weighted DAGs. 
By this we mean that for a DAG \(\dag=([n],E)\) with weights \(w\in\RR^E\) we fix a labeling 
compatible with the topological ordering induced by \(E\). This choice results in an embedding 
\(\RR^E\to\TT^{n\times n}\) for which \(c_{ij} = w(j\to i)\) and \(c_{ii} = 0\).

Following \citeauthor{Joswig:ETC}~\cite[Section 6.5.]{Joswig:ETC} we can generalize the original definition from \citeauthorandcite{Joswig.Kulas:2010}
and call a subset \(P\subseteq\TTAA^{d-1}\) a \emph{polytrope} if \(P = \tconv(C)\) for some finite 
\(C\subset\TT\PP^{d-1}\) and \(P\) is classically convex under the identification 
\(\TTAA^{d-1}\cong\mathbb{R}^{d-1}\).

A polytrope is also an ordinary polyhedron, which means it has a description 
as an intersection of classical halfspaces in the following way.
For a matrix \(C\in\mathbb{T}^{n\times n}\) its \emph{weighted digraph polyhedron} \(\Q(C)\) 
is given by the (ordinarily convex) set 
\begin{equation}
    \Q(C) = \{\, 
        x\in\TTAA^{n-1} \mid x_i - x_j \leq c_{ij} \text{ for all } 1\leq i,j\leq n, i\neq j
    \,\}.\label{eq:weighted_digraph_polyhedron}
\end{equation}

If a matrix \(C\in\TT^{n\times n}\) is the weighted adjacency matrix of a digraph \(\dag\),
calculating the shortest paths between all pairs of nodes corresponds to the evaluation of the 
polynomial \[
    C^\star = I \oplus C \oplus C^{\odot 2} \oplus \dots \oplus C^{\odot (n-1)} \oplus \dots.
\] This infinite sum converges to the expression \(C^\star = (I \oplus C)^{\odot (n-1)}\) 
if the corresponding weighted digraph has no negative cycles \cite[Lemma 3.23]{Joswig:ETC}. 
For our purposes, we focus on DAGs which means there are no cycles at all and thus \(C^\star\) 
always converges.

The matrix \(C^\star\) is called the (min-plus) \emph{Kleene star of \(C\)}; we denote by
\(\dag^\star = ([n],E^\star)\) the graph on which this matrix is supported, and by \(w^\star\in\RR^{E^\star}\) the corresponding vector of weights 
on \(\dag^\star\). In the literature, \(\dag^\star\) is sometimes called the \emph{transitive closure} of \(\dag\).

If the Kleene star \(C^\star\) converges, we have \(c_{ii} = 0\) and its entries satisfy the triangle inequality 
\(c_{ij}\leq c_{ik} + c_{kj}\). This also implies that taking the Kleene star is idempotent, \ie{} \((C^\star)^\star = C^\star\),
and multiplication by a Kleene star is also idempotent, so \(C^\star \odot C^\star = C^\star\).

\begin{proposition}[{\cite[Lemma 6]{De-La-Puente:2013}}]
    \label{lem:weighted-kleene-star-polyhedron}
    For \(C\in\mathbb{T}^{n\times n}\), we have \(\Q(C) = \Q(C^\star)\).
\end{proposition}

\begin{example}[{\cite[Example 3.\,28.]{Joswig:ETC}}]\label{ex:weighted_digraph_polyhedron}
    Consider the matrix \[
        C = \begin{pmatrix}
            1 & 4 & 0\\
            -1 & 0 & -3\\
            5 & \infty & \infty
        \end{pmatrix}
    \] with Kleene star \[
        C^\star = \begin{pmatrix}
            0 & 4 & 0\\
            -1 & 0 & -3\\
            5 & 9 & 0
        \end{pmatrix}.
    \] 

    \begin{figure}[t]
        \centering
        \begin{tikzpicture}[
        line width=1pt,
        roundnode/.style={circle, draw=black, fill=black, thin, inner sep=1pt}
    ]
    \begin{axis}[
            axis equal image,grid=both,
            ymin=-1,ymax=6,xmax=0,xmin=-5,
            xtick distance=1,
            ytick distance=1,
            enlargelimits=0.05,
            xticklabel={\color{gray}\pgfmathprintnumber\tick},
            yticklabel={\color{gray}\pgfmathprintnumber\tick},
            axis line style={draw=none}
        ]

        \node[roundnode] (a1) at (axis cs:-1,5) {};
        \node[roundnode] (a2) at (axis cs:-4,5) {};
        \node[roundnode] (a3) at (axis cs:-3,0) {};

        \node[anchor=south west,scale=1] at (a1) {$c^{\star(1)}$};
        \node[anchor=south east,scale=1] at (a2) {$c^{\star(2)}$};
        \node[anchor=north west,scale=1] at (a3) {$c^{\star(3)}$};

        \coordinate (b1) at (axis cs:-4,0);
        \coordinate (b2) at (axis cs:-1,2);

        \draw[black,thick,fill=blue,opacity=0.3] (a1.center) -- (a2.center) -- (b1.center) -- (a3.center) -- (b2.center) -- cycle;

        \draw[orange!50, thick] (axis cs:-6,3) -- (axis cs:-2,7);
    \end{axis}
\end{tikzpicture}
        \caption{The polytrope from \Cref{ex:weighted_digraph_polyhedron} together with the non-facet defining
        line defined by \(x_3 - x_2 = 9\).}
        \label{fig:ex:weighted_digraph_polyhedron}
    \end{figure}

    Evaluating \Cref{eq:weighted_digraph_polyhedron} for this particular matrix \(C^\star\) yields the 
    six defining inequalities \begin{align*}
        x_2 - x_1 &\leq -1,& x_1 - x_2 &\leq 4,& x_1 - x_3 &\leq 0, \\
        x_3 - x_1 &\leq 5,& x_3 - x_2 &\leq 9,& x_2 - x_3 &\leq -3.
    \end{align*}
    Note how the \((3,2)\)-th entry is the only off-diagonal entry not agreeing between \(C\) and \(C^\star\) 
    and how the inequality \(x_3 - x_2 \leq 9\) is redundant. This is shown as the line in \Cref{fig:ex:weighted_digraph_polyhedron}.
\end{example}

\begin{theorem}
    \label{thm:polytrope_is_weighted_digraph_polyhedron}
    A tropical polytope \(P\subset\TTAA^{n-1}\) is a polytrope if and only if 
    \[P = \tconv(C^\star) = \Q(C)\] for some \(C\in\mathbb{T}^{n\times n}\). 

    \begin{proof}
        See Section 3.\,4.\ from \citeauthorandcite{Joswig.Loho:2016}.
    \end{proof}
\end{theorem}

As a consequence of this result, we may assume that a polytrope is generated by a square matrix \(C\) such that \(c_{ii}=0\).
This also means that the construction of the weighted digraph polyhedron does not distinguish between 
a matrix \(C\) and its Kleene star $C^\star$, as seen in \Cref{lem:weighted-kleene-star-polyhedron} and \Cref{ex:weighted_digraph_polyhedron}.
On the other hand, the tropical convex hulls of a matrix and its Kleene star do differ in general. 
This can also be seen with the tropical (non-)equivalence of \(C\) and \(C^\star\) as in the following example.

\begin{example}\label{ex:degenerate-triangle}
\begin{figure}[th]
        \centering
        \begin{subfigure}[b]{.32\linewidth}
            \centering
            \begin{tikzpicture}[
                scale=1.5,
                roundnode/.style={circle, draw=black!85, fill=black!85, thin, scale=.5},
                pseudonode/.style={circle, draw=black!85, fill=black!35, scale=.5}
            ]
                \node at (.5,.5) {};
            
                \coordinate (v1) at (-1,1);
                \coordinate (v2) at (-1,0);
            
                \coordinate (v3) at (-2,-1);
                \coordinate (v4) at (-1,-1.5);
            
                \draw[black!15, fill=black!15] (v1) -- (v2) -- (v3) --++ (0,-.5) --++ (1,0) -- cycle;
            
                \draw[black!85!black, line width=2pt, loosely dotted] (v2) ++ (.5,.5) -- (v2); 
                \draw[black!85!black, line width=2pt] (v2) -- (v3); 
                \draw[black!85!black, line width=2pt, loosely dotted] (v1) ++ (-1,0) -- (v1);
                \draw[black!85!black, line width=2pt] (v1) -- (v4);
                \draw[black!85!black, line width=2pt, loosely dotted] (v1) ++ (.5,.5) -- (v1);
            
                \node[pseudonode] at (v2) {};
                \node[roundnode] at (v1) {};
            
                \node[black!85!black, above=of v1, anchor=north] {$c^{(1)}$};
            \end{tikzpicture}
            \caption{\(\tconv(C)\)}
        \end{subfigure}%
        \begin{subfigure}[b]{.32\linewidth}
            \centering
            \begin{tikzpicture}[
                scale=1.5,
                roundnode/.style={circle, draw=black!85, fill=black!85, thin, scale=.5},
                pseudonode/.style={circle, draw=black!85, fill=black!35, scale=.5}
            ]
                \node at (.5,.5) {};
            
                \coordinate (v1) at (-1,1);
                \coordinate (v2) at (-1,0);
            
                \coordinate (v3) at (-2,-1);
                \coordinate (v4) at (-1,-1.5);
            
                \draw[black!15, fill=black!15] (v1) -- (v2) -- (v3) --++ (0,-.5) --++ (1,0) -- cycle;
            
                \draw[black!85!black, line width=2pt, loosely dotted] (v2) ++ (.5,.5) -- (v2); 
                \draw[black!85!black, line width=2pt] (v2) -- (v3); 
                \draw[black!85!black, line width=2pt, loosely dotted] (v1) ++ (-1,0) -- (v1) --++(.5,0);
                \draw[black!85!black, line width=2pt, loosely dotted] (v1) ++ (0,.5) -- (v4);
                \draw[black!85!black, line width=2pt] (v2) -- (v4);
            
                \node[pseudonode] at (v2) {};
            \end{tikzpicture}
            \caption{\(\Q(C)\)}
        \end{subfigure}%
        \begin{subfigure}[b]{.32\linewidth}
            \centering
            \begin{tikzpicture}[
                scale=1.5,
                roundnode/.style={circle, draw=black!85, fill=black!85, thin, scale=.5},
                pseudonode/.style={circle, draw=black!85, fill=black!35, scale=.5}
            ]
                \node at (.5,.5) {};
            
                \coordinate (v1) at (-1,0);
                \coordinate (v2) at (-1,0);
            
                \coordinate (v3) at (-2,-1);
                \coordinate (v4) at (-1,-1.5);
            
                \draw[black!15, fill=black!15] (v1) -- (v2) -- (v3) --++ (0,-.5) --++ (1,0) -- cycle;
            
                \draw[black!85!black, line width=2pt, loosely dotted] (v2) ++ (.5,.5) -- (v2); 
                \draw[black!85!black, line width=2pt] (v2) -- (v3); 
                \draw[black!85!black, line width=2pt, loosely dotted] (v1) ++ (-1,0) -- (v1);
                \draw[black!85!black, line width=2pt] (v1) -- (v4);
                \draw[black!85!black, line width=2pt, loosely dotted] (v1) ++ (.5,.5) -- (v1);
            
                \node[pseudonode] at (v2) {};
                \node[roundnode] at (v1) {};
     
                \node[black!85!black, above=of v1, anchor=north] {$c^{\star(1)}$};
            \end{tikzpicture}
            \caption{\(\Q(C^\star) = \tconv(C^\star)\)}
        \end{subfigure}
        \caption{Comparison of the hyperplane arrangements in \Cref{ex:degenerate-triangle}.}
        \label{fig:icecream-cone-arrangement}
    \end{figure}  
    Continuing from \Cref{ex:affine-covectors}, we consider the case \(c > a + b\) which yields the Kleene star \[
        C^\star = \begin{pmatrix}
            0 & \infty & \infty \\
            a & 0 & \infty \\
            a+b & b & 0
        \end{pmatrix}.
    \] The corresponding hyperplane arrangements for \(\mathcal{T}(C)\) and \(\mathcal{T}(C^\star)\) are shown in
    \Cref{fig:icecream-cone-arrangement}. The point configurations given by \(C\) and \(C^\star\) are not tropically 
    equivalent since their tropical combinatorial types \(\tc(C)\) and \(\tc(C^\star)\) are not isomorphic as posets.

    We can see in \Cref{fig:tc-poset} that \(\tc(C)\) has five atoms (i.e., elements covering the mininum)
    corresponding to the five chambers. 
    On the other hand, \(\tc(C^\star)\) has four atoms corresponding to the four chambers
    in \hyperlink{fig:icecream-cone-arrangement}{\Cref*{fig:icecream-cone-arrangement}.(c)}.
    Compared to \(\tc(C)\), the covectors \((\varnothing,1,23)\) and \((1,1,23)\) are missing
    and \((\varnothing,1,123)\) gets identified with \((\varnothing,12,23)\).
    
    Yet, \(\Q(C)\) does occur in the hyperplane arrangement of \(\mathcal{T}(C)\) and 
    as a subpolytope of \(\tconv(C)\), namely as the two-dimensional unbounded cell in \Cref{fig:icecream-cone-arrangement}.
    \begin{figure}[bh]
       \centering 
%
%
%
%
%
%
%
%
\newsavebox\vuuu
\begin{lrbox}{\vuuu}
\begin{tikzcd}[font = \Large, blue]
    1 \arrow[d] & 2 \arrow[d] & 3 \arrow[d] \\
    1 & 2 & 3
\end{tikzcd}
\end{lrbox}
\newsavebox\vzdu
\begin{lrbox}{\vzdu}
\begin{tikzcd}[font = \Large,black!60]
    1 \arrow[dr]& 2 \arrow[d] & 3 \arrow[d] \\
    1 & 2 & 3
\end{tikzcd}
\end{lrbox}
\newsavebox\vzud
\begin{lrbox}{\vzud}
\begin{tikzcd}[font = \Large,black!60]
    1 \arrow[dr] & 2 \arrow[dr] & 3 \arrow[d] \\
    1 & 2 & 3
\end{tikzcd}
\end{lrbox}
\newsavebox\vuzd
\begin{lrbox}{\vuzd}
\begin{tikzcd}[font = \Large,black!60]
    1 \arrow[d] & 2 \arrow[dr] & 3 \arrow[d] \\
    1 & 2 & 3
\end{tikzcd}
\end{lrbox}
\newsavebox\vzzt
\begin{lrbox}{\vzzt}
\begin{tikzcd}[font = \Large,black!60]
    1 \arrow[drr] & 2 \arrow[dr] & 3 \arrow[d] \\
    1 & 2 & 3
\end{tikzcd}
\end{lrbox}
\newsavebox\vuudd
\begin{lrbox}{\vuudd}
\begin{tikzcd}[font = \Large,blue]
    1 \arrow[d] & 2 \arrow[d]\arrow[dr] & 3 \arrow[d] \\
    1 & 2 & 3
\end{tikzcd}
\end{lrbox}
\newsavebox\vuud
\begin{lrbox}{\vuud}
\begin{tikzcd}[font = \Large]
    1 \arrow[d]\arrow[dr] & 2 \arrow[dr] & 3 \arrow[d] \\
    1 & 2 & 3
\end{tikzcd}
\end{lrbox}
\newsavebox\vudu
\begin{lrbox}{\vudu}
\begin{tikzcd}[font = \Large,blue]
    1 \arrow[d]\arrow[dr] & 2 \arrow[d] & 3 \arrow[d] \\
    1 & 2 & 3
\end{tikzcd}
\end{lrbox}
\newsavebox\vzut
\begin{lrbox}{\vzut}
\begin{tikzcd}[font = \Large,black!60]
    1 \arrow[drr]\arrow[dr] & 2 \arrow[dr] & 3 \arrow[d] \\
    1 & 2 & 3
\end{tikzcd}
\end{lrbox}
\newsavebox\vzdd
\begin{lrbox}{\vzdd}
\begin{tikzcd}[font = \Large,black!60]
    1 \arrow[dr]& 2 \arrow[dr]\arrow[d] & 3 \arrow[d] \\
    1 & 2 & 3
\end{tikzcd}
\end{lrbox}
\newsavebox\vuzt
\begin{lrbox}{\vuzt}
\begin{tikzcd}[font = \Large,black!60]
    1 \arrow[d]\arrow[drr] & 2 \arrow[dr]\arrow[d] & 3 \arrow[d] \\
    1 & 2 & 3
\end{tikzcd}
\end{lrbox}
\newsavebox\vuut
\begin{lrbox}{\vuut}
\begin{tikzcd}[font = \Large]
    1 \arrow[d]\arrow[drr] & 2 \arrow[dr] & 3 \arrow[d] \\
    1 & 2 & 3
\end{tikzcd}
\end{lrbox}
\newsavebox\vudd
\begin{lrbox}{\vudd}
\begin{tikzcd}[font = \Large, blue]
    1 \arrow[d]\arrow[dr] & 2 \arrow[dr]\arrow[d] & 3 \arrow[d] \\
    1 & 2 & 3
\end{tikzcd}
\end{lrbox}

\begin{tikzpicture}[x  = {(6em, 0em)},
                    y  = {(0em, 4em)},
                    scale = 1,
                    color = {lightgray}]
  \coordinate (vbottom) at (0, 0);
  
  \coordinate (v111)   at (   0, 1);
  \coordinate (v021)   at (   2, 1);
  \coordinate (v012)   at (   1, 1);
  \coordinate (v102)   at (  -1, 1);
  \coordinate (v003)   at (  -2, 1);
  
  \coordinate (v121)   at ( 1.5, 2);
  \coordinate (v022)   at ( 2.5, 2);
  \coordinate (v013)   at (- .5, 2);
  \coordinate (v112)   at (  .5, 2);
  \coordinate (v112_2) at (-1.5, 2);
  \coordinate (v103)   at (-2.5, 2);
  
  \coordinate (v122)   at (   1, 3);
  \coordinate (v113)   at (  -1, 3);
  
  \coordinate (vtop) at (0, 4);

  \definecolor{vertexcolor_bottom}{rgb}{ 1 1 1 }
  \definecolor{vertexcolor_proper}{rgb}{ 1 1 1 }
  \definecolor{vertexcolor_top}{rgb}{ 0 0 0 }

  \colorlet{vertexbordercolor}{black}
  \colorlet{vertexcolor_nonface}{black!50}
  \colorlet{vertexcolor_tconv}{black}
  \colorlet{vertexcolor_polytrope}{blue}

  \tikzstyle{vertexstyle_bottom} = [text=black, inner sep=2pt, rectangle, rounded corners=3pt,fill=vertexcolor_bottom, draw=vertexbordercolor,]
  
  \tikzstyle{vertexstyle_111} = [text=vertexcolor_polytrope, inner sep=2pt, rectangle, rounded corners=3pt,fill=vertexcolor_proper, draw=vertexcolor_polytrope,]
  \tikzstyle{vertexstyle_021} = [text=vertexcolor_nonface, inner sep=2pt, rectangle, rounded corners=3pt,fill=vertexcolor_proper, draw=vertexcolor_nonface,]
  \tikzstyle{vertexstyle_012} = [text=vertexcolor_nonface, inner sep=2pt, rectangle, rounded corners=3pt,fill=vertexcolor_proper, draw=vertexcolor_nonface,]
  \tikzstyle{vertexstyle_102} = [text=vertexcolor_nonface, inner sep=2pt, rectangle, rounded corners=3pt,fill=vertexcolor_proper, draw=vertexcolor_nonface,]
  \tikzstyle{vertexstyle_003} = [text=vertexcolor_nonface, inner sep=2pt, rectangle, rounded corners=3pt,fill=vertexcolor_proper, draw=vertexcolor_nonface,]
  
  \tikzstyle{vertexstyle_121} = [text=vertexcolor_polytrope, inner sep=2pt, rectangle, rounded corners=3pt,fill=vertexcolor_proper, draw=vertexcolor_polytrope,]
  \tikzstyle{vertexstyle_022} = [text=vertexcolor_nonface, inner sep=2pt, rectangle, rounded corners=3pt,fill=vertexcolor_proper, draw=vertexcolor_nonface,]
  \tikzstyle{vertexstyle_013} = [text=vertexcolor_nonface, inner sep=2pt, rectangle, rounded corners=3pt,fill=vertexcolor_proper, draw=vertexcolor_nonface,]
  \tikzstyle{vertexstyle_112} = [text=vertexcolor_tconv, inner sep=2pt, rectangle, rounded corners=3pt,fill=vertexcolor_proper, draw=vertexcolor_tconv,]
  \tikzstyle{vertexstyle_112_2} = [text=vertexcolor_polytrope, inner sep=2pt, rectangle, rounded corners=3pt,fill=vertexcolor_proper, draw=vertexcolor_polytrope,]
  \tikzstyle{vertexstyle_103} = [text=vertexcolor_nonface, inner sep=2pt, rectangle, rounded corners=3pt,fill=vertexcolor_proper, draw=vertexcolor_nonface,]
  
  \tikzstyle{vertexstyle_122} = [text=vertexcolor_polytrope, inner sep=2pt, rectangle, rounded corners=3pt,fill=vertexcolor_proper, draw=vertexcolor_polytrope,]
  \tikzstyle{vertexstyle_113} = [text=vertexcolor_tconv, inner sep=2pt, rectangle, rounded corners=3pt,fill=vertexcolor_proper, draw=vertexcolor_tconv,]

  \tikzstyle{vertexstyle_top} = [text=black, inner sep=2pt, rectangle, rounded corners=3pt,fill=vertexcolor_top, draw=vertexbordercolor,]

  \tikzstyle{edgestyle_nonface} = [thick,color=black!50]
  \tikzstyle{edgestyle_tconv} = [thick,color=black]
  \tikzstyle{edgestyle_polytrope} = [thick,color=blue]


  \foreach \i/\k/\style in {
    bottom/111/nonface,
    bottom/012/nonface,
    bottom/021/nonface,
    bottom/102/nonface, 
    bottom/003/nonface,
    121/111/polytrope,
    121/021/nonface,
    013/012/nonface, 
    013/003/nonface,
    022/012/nonface,
    022/021/nonface,
    112/102/nonface,
    112/012/nonface,
    112_2/111/polytrope,
    112_2/102/nonface,
    103/003/nonface,
    103/102/nonface,
    122/121/polytrope,
    122/022/nonface,
    122/112/tconv,
    122/112_2/polytrope,
    113/103/nonface,
    113/013/nonface,
    113/112/tconv,
    122/top/nonface,
    113/top/nonface%
  } {
   \draw[edgestyle_\style] (v\i) -- (v\k);
  }

  \foreach \i/\label in {
  	bottom/ ,
	top/ ,
	111/{\usebox\vuuu},
	021/{\usebox\vzdu},
        012/{\usebox\vzud},
        102/{\usebox\vuzd},
        003/{\usebox\vzzt},
        121/{\usebox\vudu},
        022/{\usebox\vzdd},
        013/{\usebox\vzut},
        112/{\usebox\vuud},
        112_2/{\usebox\vuudd},
        103/{\usebox\vuzt},
        122/{\usebox\vudd},
        113/{\usebox\vuut}
  } {
    \node at (v\i) [vertexstyle_\i] {\adjustbox{scale=.4}{\label}};
  }

\end{tikzpicture}
       \caption{Hasse diagram of \(\tc(C)\) from \Cref{ex:degenerate-triangle}. The subposet corresponding 
       to \(\tconv(C)\) is marked in black and blue, while the subposet of \(\Q(C)\) is marked in blue.}
       \label{fig:tc-poset}
     \end{figure}
\end{example}

Note that a weighted digraph polyhedron can be expressed as the set of points 
\({x\in\TTPP^{d-1}}\) satisfying the tropical linear inequalities \(x\leq C\odot x\) since \[
    x_i - x_j \leq c_{ij} \iff x_i \leq x_j + c_{ij}
\] for all \(i\in[n]\) and since the right inequality holds for all \(j\in[n]\),  
\begin{equation}
    x_i\leq\min_j \left\{\,x_j + c_{ij}\,\right\} = c_i\odot x \label{eq:tropical-inequalities}
\end{equation} also holds.

Everything until this point has set the stage for our first main result, which now follows without much 
effort. This allows us to recover the polytrope for a given point configuration from 
the combinatorial data of the tropical hyperplane arrangement. Recall that a \emph{perfect matching} in a graph is a subset of the edges that covers every node exactly once.
\begin{theorem}\label{thm:polytrope-cell}
    Let \(C\in\TT^{n\times n}\) be supported on a DAG \(\dag\). 
    The polytrope \(\tconv(C^\star)\subseteq\tconv(C)\) coincides with the closed cell \(\sector_L(C)\) 
    such that \(L\) is the perfect matching \((1,2,\dots,n)\) in the support graph 
    \(\supportgraph[C]\).

    \begin{proof}
        By assumption, \(\supportgraph[C]\) contains the perfect matching \((1,2,\dots,n)\)
        as \(C\) is assumed to be supported on a DAG. Also, \(\sector_L(C)\subseteq\tconv(C)\) 
        by \Cref{thm:covdec-tconv-cells} which means we are left to prove that 
        \(\sector_L(C)=\tconv(C^\star)\).
        
        Since \(C\) is supported on a DAG, the polytrope \(P=\tconv(C^\star)\) is non-empty
        and can be described as the set of points satisfying 
        the inequalities in \eqref{eq:tropical-inequalities} for all \(i,j\in[n]\). 
        Rearranging those inequalities gives \[
        	x_j \geq \max_{1\leq i\leq n}\{-c_{ij} + x_i\} = \alpha_j(x)
         \] for every \(j\in[n]\).
        Thus, any point \(x\in P\) is inside every closed sector \(\sector_j(v^{(j)})\) for each \(j\in[n]\).
        By \eqref{eq:cap-sectors-cup-covector}, the intersection of all these sectors 
        is the covector with edges \((j,j)\) for every \(j\in[n]\).
    \end{proof}
\end{theorem}

Another approach was taken by \citeauthorandcite{Joswig.Schröter:2019}.  
Instead, they characterized combinatorial types of polytropes in
terms of \emph{central triangulations} of the following lower-dimensional point configuration.

\begin{definition}
  For a digraph \(\dag\), the associated \emph{fundamental polytope} is defined as 
  \[
    \fundamentalpolytope[\dag] 
    = \conv\left(
        \left\{\, e_i - e_j \mid j\to i\in\dag \,\right\} \cup \{0\}
    \right).
  \]
\end{definition}

The fundamental polytope \(\fundamentalpolytope\) provides an appropriate restriction of the regular subdivisions
of \(\Newton(\defpoly[C])\) when \(C\in\TT^{n\times n}\) is supported on a DAG \(\dag\). 
The relationship between \(\fundamentalpolytope\) and \(\Newton(\defpoly[C])\) is of the following nature.
Recall from \Cref{rem:sectors-monomials} the correspondence between the open covector cells of \(\CovDec(C)\)
and the monomials of \(\defpoly[C]\). For a covector \(L\) with coarse type \(a\) such that \(a_j = 0\), \(a_i = 2\)
and otherwise \(a_k = 1\), \(L\) must be (since \(\dag\) contains no cycles) of the form \[
    \begin{tikzcd}[font = \small, row sep=1.75em, column sep=2em]
        1 \arrow[d] & 2 \arrow[d] & \dots & j\arrow[dr] & i\arrow[d] & \dots & n \arrow[d]\\
        1 & 2 & \dots & j & i & \dots & n
    \end{tikzcd}.
\] The coefficient of \(x^a\) in \(\defpoly[C]\) is given by the weight of the above covector \(L\) and is precisely \(c_{ij}\),
which means that the vertex \(a\) is also lifted to height \(c_{ij}\) in the regular subdivision of \(\Newton(\defpoly[C])\).
But then, \(a = e_i - e_j + \mathbf{1}\), which gives \(\fundamentalpolytope\) as a translated subpolytope
of \(\Newton(\defpoly[C])\) with agreeing height function.

Motivated by this correspondence, we say that a subdivision of the fundamental polytope is \emph{central} if
the origin \(\mathbf{0}\) is a vertex of every maximal cell. This terminology comes from the fact that
for the complete digraph \(K_n\), the origin is the unique interior lattice point of both \(\fundamentalpolytope\)
and \(\Newton(\defpoly)\).

\begin{example}
\begin{figure}[th]
    \centering
    \begin{tikzpicture}[
        roundnode/.style={circle, draw=black!85, fill=black!85, thin, scale=0.5},
        pseudonode/.style={circle, draw=black!85, fill=white, scale=0.5}
    ]
        \coordinate (0) at (0,0);
        \coordinate (12) at (-1,1,0);
        \coordinate (21) at (1,-1,0);
        \coordinate (13) at (-1,0,1);
        \coordinate (31) at (1,0,-1);
        \coordinate (23) at (0,-1,1);
        \coordinate (32) at (0,1,-1);
        
        \draw[black, fill=black!15] (0) -- (12) -- (13) -- (23) -- cycle;
        \draw[black] (0) -- (13);
    
        \node[roundnode] at (0) {};
        \node[roundnode] at (12) {}; \node[anchor=south east] at (12) {$e_2 - e_1$};
        \node[roundnode] at (21) {};
        \node[roundnode] at (13) {}; \node[anchor=east] at (13) {$e_3 - e_1$};
        \node[roundnode] at (31) {}; 
        \node[roundnode] at (23) {}; \node[anchor=north] at (23) {$e_3 - e_2$};
        \node[roundnode] at (32) {};
    \end{tikzpicture}
    \hspace{2em}
    \begin{tikzpicture}[
        roundnode/.style={circle, draw=black!85, fill=black!85, thin, scale=0.5},
        pseudonode/.style={circle, draw=black!85, fill=white, scale=0.5}
    ]
        \coordinate (0) at (0,0);
        \coordinate (12) at (-1,1,0);
        \coordinate (21) at (1,-1,0);
        \coordinate (13) at (-1,0,1);
        \coordinate (31) at (1,0,-1);
        \coordinate (23) at (0,-1,1);
        \coordinate (32) at (0,1,-1);
        
        \draw[black, fill=black!15] (0) -- (12) -- (13) -- (23) -- cycle;
        \draw[black] (12) -- (23);
    
        \node[roundnode] at (0) {};
        \node[roundnode] at (12) {}; \node[anchor=south east] at (12) {$e_2 - e_1$};
        \node[roundnode] at (21) {};
        \node[roundnode] at (13) {}; \node[anchor=east] at (13) {$e_3 - e_1$};
        \node[roundnode] at (31) {}; 
        \node[roundnode] at (23) {}; \node[anchor=north] at (23) {$e_3 - e_2$};
        \node[roundnode] at (32) {};
    \end{tikzpicture}
    \caption{Two triangulations of the fundamental polytope \(\fundamentalpolytope\) for \(G = \kappa_3\).
    The left triangulation is central, while the right triangulation is not.}
    \label{fig:example-central-triangulation}
    \end{figure}

    \Cref{fig:example-central-triangulation} shows two subdivisions for \(\fundamentalpolytope\) where \(G=\kappa_3\).
    The first triangulation is central because the origin is a vertex of every maximal cell. This triangulation appears 
    for all weights \(w\in\RR^E\) such that \(w(1\to 3) < w(1\to 2) + w(2\to 3)\).

    On the other hand, if \(w(1\to 3) > w(1\to 2) + w(2\to 3)\) we obtain the second triangulation which is not central.
\end{example}

\begin{remark}\label{rem:root-subsets}
    The fundamental polytope \(\fundamentalpolytope[\dag]\) corresponds to a subset of 
    the \(A_n\) root system together with the origin. In that case, a matrix \(C\in\TT^{n\times n}\) 
    supported on a DAG \(\dag\) induces a regular subdivision on the \(A_n\) root system which involves 
    precisely the vertices of \(\fundamentalpolytope\) while omitting the other vertices by lifting 
    them to infinite height.
\end{remark}

Recall that a directed graph \(\dag\) is \emph{transitive} if 
for any pair of edges \(j\to k\) and \(k\to i\), \(\dag\) also contains the edge \(j\to i\).
That is, viewing the edges of \(\dag\) as a relation on the nodes gives a transitive relation.
Sometimes, such a graph is also called \emph{transitively closed}.
We now adapt the proof of \cite[Theorem 22]{Joswig.Schröter:2019} to the case
of arbitrary support.
\begin{theorem}\label{thm:central-triangulations}
    The tropical combinatorial types of full-dimensional polytropes in \(\TTAA^{n-1}\) are in
    bijection with the regular central subdivisions of \(\fundamentalpolytope\) where \(G\)
    ranges over transitive directed graphs \(\dag\) on \(n\) nodes.

    \begin{proof}
        By \Cref{thm:polytrope_is_weighted_digraph_polyhedron}, we know that for any given 
        non-empty polytrope \(P\subseteq\TTAA^{d-1}\) there is a matrix \(C\in\TT^{n\times n}\) with 
        \(P = \tconv(C)\) and \(C = C^\star\). Denote by \(\dag\) the underlying directed graph of \(C\).
        By \Cref{prop:mixed-regular-subdivision-newton}, \(P\) is dual to the regular subdivision \(\subdivision\) of 
        \(\Newton(\defpoly[C])\subseteq n\cdot\Delta_{n-1}\) induced by the coefficients of \(\defpoly[C]\). 
        Thus, the height \(h(a)\) of a vertex \(a = \sum_{i=1}^n a_i e_i \in \Newton(\defpoly[C])\) is given by 
        \begin{equation}\label{eqn:height-function}
            h(a) = \min_L \sum_{j = 1}^n c_{ij}
        \end{equation}
        where \(L\) is a covector corresponding to a full-dimensional cell in \(\CovDec(C)\) such that 
        \(\#L_i = a_i\). In particular, any such covector \(L\) has a unique edge \(j\to i\) for each \(j\in[n]\).

        Any Kleene star matrix satisfies the triangle inequality \(c_{ij}\leq c_{ik} + c_{kj}\)
        and has zero diagonal.
        This means the central vertex \(\mathbf{1}\) is lifted to height zero since the
        weight of the matching \((1,2,\dots,n)\) evaluates to the sum of the diagonal entries of \(C\) in \eqref{eqn:height-function}.
        
        If there was another matching \(L'\) with lower weight, then \(C\) would have a negative cycle,
        which contradicts \(P\) being non-empty. Using the triangle inequality \(0 = c_{ii}\leq c_{ij} + c_{ji}\),
        we can see that the central vertex lies below each line spanned by any two vertices 
        \(a,a'\in\Newton(\defpoly[C])\). Since \(\mathbf{1}\) is the unique interior lattice point of 
        \(\Newton(\defpoly[C])\) with respect to \(n\cdot\Delta_{n-1}\), this means \(\mathbf{1}\)
        is a vertex of every maximal cell. Thus, \(\subdivision\) is a 
        central mixed regular subdivision of \(\Newton(\defpoly[C])\).

        Now, the matrix \(C\) is the weighted adjacency matrix of a graph \(\dag\).
        Since \(\fundamentalpolytope\) can be seen as a translated subpolytope of \(\Newton(\defpoly[C])\), we can restrict \(S\) to \(\fundamentalpolytope\). 
        This is precisely the regular subdivision of \(\fundamentalpolytope\) induced by \(C\) 
        since the height function on \(\fundamentalpolytope\) agrees with the height function on \(\Newton(\defpoly[C])\).

        This means the maximal-dimensional cells of \(\Newton(\defpoly[C])\) are in bijection with 
        the maximal dimensional cells of \(\fundamentalpolytope\). As a consequence,
        the central mixed subdivisions of \(\Newton(\defpoly[C])\) are mapped bijectively 
        to the central subdivisions of \(\fundamentalpolytope\).
    \end{proof}
\end{theorem}

Connecting back to \Cref{rem:root-subsets}, we can see that the combinatorial types of polytropes
correspond to regular central subdivisions of those subsets of the \(A_n\) root system corresponding to 
transitive DAGs on \(n\) nodes. 

\section{Moduli space of polytropes from DAGs}\label{sec:coarse-classification}
We have established two characterizations for tropical equivalence of polytropes, via special subcomplexes
in regular subdivisions of root polytopes and duals of central regular subdivisions of fundamental polytopes.
This suggests that certain subfans of the secondary fans \(\sfan(\rootpolytope)\) resp.\ 
\(\sfan(\fundamentalpolytope)\) classify combinatorial types of polytropes.

If \(\dag\) is transitive, this subfan is cut out by the region of Kleene star matrices on \(\dag\),
which in the case of the complete graph has been noted by \citeauthorandcite{Tran:2016}.
If \(\dag\) is not transitive, this region is not useful, as taking the Kleene star of a matrix 
inserts edges in the supporting graph which enlarges both \(\rootpolytope[\dag^\star]\)
and \(\fundamentalpolytope[\dag^\star]\) as compared to \(\rootpolytope\) and \(\fundamentalpolytope\).
In other words, there are no Kleene star matrices supported precisely on \(\dag\).

At the same time, specific choices of edge weights in \(C\) might give a polytrope \(\Q(C)\) exhibiting 
the same face structure as a polytrope supported on a subgraph \(\subdag\subset\dag\).
We can illustrate this phenomenon already with our running example.
\begin{example}\label{ex:redundant-face}
    There are two possible cases for the directed triangle \(\kappa_3\) with weight matrix \[
        C = \begin{pmatrix}
            0 & \infty & \infty \\
            a & 0 & \infty \\
            c & b & 0
        \end{pmatrix}.
    \]
        
    If \(c \geq a+b\), then one can see in \Cref{fig:icecream-cone-arrangement} that the hyperplane 
    \(x_3 - x_2 = c\) does not define a facet and is redundant. Thus, one can replace the value of \(c\)
    by \(\infty\) and delete the edge \(1\to3\) while keeping the resulting polytrope intact.

    On the other hand, in \Cref{fig:hyperplane-arrangement-sectors}, the case \(c < a+b\) is shown.
    There, \(x_3 - x_2 = c\) is a facet-defining hyperplane.
\end{example}

\Cref{ex:redundant-face} suggests that an edge not contributing to a shortest path between any pair of nodes
does not contribute to the polytrope. Such a path is not \emph{critical} in the sense of
\citeauthorandcite{amendola2022conditional}.
Thus, the solution to finding a non-redundant weight matrix \(C\) 
involves identifying the graph giving the minimal non-redundant facet description.

\subsection{Transitive reduction}
The characterization of non-redundant facets of \(\Q(C)\) and edges of \(\dag\) relies on the 
interplay between \emph{transitive closure} and \emph{covering relations} of the edge set \(E\) 
seen as a binary relation.

\begin{definition}\label{def:transitive}
    For a given DAG \(\dag\), the \emph{transitive reduction} \(\dag^{\trred}\) contains an edge \(j\to i\in G\) 
    if there is no other path between \(j\) and \(i\).

    For a DAG \(\dag\) with weights \(w\in\RR^E\), the \emph{weighted transitive reduction} \(\dag^\flat_w\) 
    is the digraph containing an edge \(j\to i\) with weight \(w(j\to i)\) whenever 
    \(j\to i\) is the unique shortest path between \(j\) and \(i\) in \(\dag\).
    
    We denote by \(w^\flat = w\vert_{\dag^\flat_w}\) the vector of weights on \(\dag^\flat_w\).
    If \(C\) is the weighted adjacency matrix of a DAG \(\dag\), denote by \(C^\flat\) the weighted adjacency matrix of \(\dag^\flat_w\).
\end{definition}

The weighted transitive reduction \(G^\flat_w\) has appeared in \cite{Gissibl:PhD-2018} as the \emph{minimum max-linear DAG}.
The sets of edges \(E^{\trred}\) of \(\dag^{\trred}\) \resp \(E^\flat_w\) of \(\dag^\flat_w\) are called the 
\emph{covering relations} \resp \emph{weighted covering relations} of \(\dag\).

\begin{remark}
    The transitive reduction \(\dag^{\trred}\) is the smallest directed graph with the same 
    reachability relation as \(\dag\).
    The weighted transitive reduction imposes the additional constraint that the appearing edges 
    already realize the shortest paths between their endpoints uniquely, thus we always have 
    \(\dag^\flat_w\supseteq \dag^{\trred}\).
\end{remark}

\begin{example}\label{ex:transitive-reduction}
    Continuing from \Cref{ex:redundant-face}, the only strict subgraph of the directed triangle~\(\dag=\kappa_3\)
    with the same reachability between all nodes is the path graph 
    \[
    	\subdag\colon\begin{tikzcd}[row sep=.5em]
            1 \arrow[r,"a"] & 2 \arrow[r,"b"] & 3
        \end{tikzcd}
    \] since for both \(\dag\) and \(\subdag\), any node \(i\) is accessible from 
    any other node \(j\) if \(i > j\). Thus, for the transitive reduction we have \(\dag^{\trred} = \subdag\).
    If \(c \geq a + b\), then \(\dag^{\trred}\) and \(\dag^\flat_{(a,b,c)}\) agree since \(1\to 2\to 3\) is a shortest path. 
    On the other hand if \(c < a + b\), then \(1\to 3\) is the unique shortest path connecting
    its endpoints, which means \(\dag = \dag^\flat_{(a,b,c)} \supsetneq \dag^{\trred}\).
\end{example}

It has been shown by \citeauthor{Aho.Garey.ea:1972}~\cite{Aho.Garey.ea:1972} that for finite unweighted DAGs 
the transitive reduction \(\dag^{\trred}\) always exists and is unique. We can extend this to finite weighted DAGs.
For this, we recall the following definition from \citeauthorandcite{Joswig.Schröter:2022}.

\begin{definition}
   If \(T^{(s)}\) is a spanning tree for the set of vertices reachable from \(s\) (including \(s\) itself),
   we say that \(T^{(s)}\) is a \emph{shortest-path tree at (source) \(s\)} in a DAG \(\dag\) with weights \(w\in\RR^E\)
   if \(s\) has no incoming edges and for every edge \(j\to i\), \[
       w(j\to i) + p_j \not< p_i,
   \] where \(p_i\) is the length of the path \(s\rightsquigarrow i\) in \(T^{(s)}\).
\end{definition}

\begin{proposition}\label{prop:weighted-transitive-reduction-trees}
  The weighted transitive reduction \(\dag^\flat_w\) of an weighted DAG \(\dag\) with weights \(w\in\RR^E\)
  always exists and is unique. It is given by \[
      \dag^\flat = \bigcup_{s\in[n]}\bigcap\{\, 
        T^{(s)} \mid T^{(s)} \text{ is a shortest-path tree at }s 
        \text{ for } w
      \,\}.
    \]
    \begin{proof}
        The right hand side \(\subdag\) contains only edges in \(\dag\) that are contained in a shortest path connecting two nodes,
        thus \(\subdag\supseteq\dag^\flat\). On the other hand, consider an edge \(j\to i\notin H\) with \(i\) reachable 
        from \(j\). This means that \(j\to i\) is not contained in every shortest-path tree at \(j\), so there exists another 
        shortest path connecting \(j\) to \(i\). But this means \(j\to i\notin\dag^\flat\), which proves the entire statement.
    \end{proof}
\end{proposition}

The following lemma gives the condition under which an edge does not contribute a facet to \(\Q(C)\)
and thus can be deleted from \(\dag\).
\begin{lemma}\label{lem:transitive-reduction-step}
    Let \(\dag\) be a DAG with weights \(w\in\RR^E\) and fix a labeling of \(\dag\).
    Then, the following are equivalent:
    \begin{enumerate}
        \item \(\Q(C) = \Q(C')\) for all \({w'\in w+\RR_{>0}e_{j\to i}}\),
        \item there exists \({w'\in w+\RR_{>0}e_{j\to i}}\) such that \(\Q(C) = \Q(C')\),
        \item \(j\to i\notin\dag^\flat\),
    \end{enumerate}
    where \(C\) and \(C'\) denote the weighted adjacency matrices arising from \(w\) resp.\ \(w'\).
	
    \begin{proof}
        We may assume \wLOG\ that \(w(j\to i) = w^\star(j\to i)\). The implication \((1)\implies(2)\) is immediate. 
        For the converse direction, if \(w'\in w+\RR_{>0}e_{j\to i}\) as required exists, 
        then the hyperplane~\({\{x_i - x_j = c_{ij}\}}\) does not define a facet of \(\Q(C)\), 
        which means that \(\Q(C) = \Q(C')\) for any value \(w'(j\to i) \geq w(j\to i)\).

        If \(w'(j\to i)\geq w(j\to i)\) but \(\Q(C) = \Q(C')\) for some \(w'\in\RR^E\), we know that 
        there exists another path with weight \(w(j\to i)\) from \(j\) to \(i\) in \(\dag\) with weights \(w'\)
        since \({(w')^\star(j\to i) = w(j\to i)}\) has to be satisfied. 
        Thus, \(j\to i\) was not the unique shortest path. Since \(\dag\) is acyclic, 
        any other path does not use the edge \(j\to i\). 
        The same reasoning works in reverse, proving the equivalence \((2)\iff(3)\).
    \end{proof}
\end{lemma}

\begin{corollary}\label{cor:transitive-reduction-step}
    Let \(\dag\) be a DAG with weights \(w\in\RR^E\), \(\dag' = \dag\setminus j\to i\) be the DAG 
    obtained from \(\dag\) by deleting the edge \(j\to i\) and denote by \(C\) the weight matrix obtained from \(w\).
    Then, \(\Q(C) = \Q(C\vert_{\dag'})\) if and only if \(j\to i\notin\dag^\flat_w\).

    \begin{proof}
        The reasoning is the same as in the proof for \Cref{lem:transitive-reduction-step}, 
        if \(j\to i\notin\dag^\flat_w\) then the \(x_i-x_j\)-hyperplane does not define a facet of \(\Q(C)\). 
        Thus removing this hyperplane does not change \(\Q(C)\) and thus \(\Q(C)=\Q(C\vert_{\dag'})\).
        If \(\Q(C) = \Q(C\vert_{\dag'})\) we may assume \wLOG\ that \(w(j\to i) = w^\star(j\to i)\) which means that 
        there is another shortest path from \(j\) to \(i\). But this means \(j\to i\notin\dag^\flat_w\).
    \end{proof}
\end{corollary}
In other words, \Cref{lem:transitive-reduction-step} describes that under the stated conditions the polytrope \(\Q(C)\) 
stays the same along rays in certain coordinate directions \(e_{j\to i}\). 
After compactifying this ray by adding a point at infinity, which corresponds to the edge \(j\to i\) being deleted from \(\dag\),
\Cref{cor:transitive-reduction-step} asserts that the polytrope does not change.

On the other hand, starting with a fixed \(\dag\) with weights \(w\), we may apply \Cref{lem:transitive-reduction-step} 
succcesively to every edge \(j\to i\in\dag\setminus\dag^\flat\) to get a pointed polyhedral cone. Taking into account the 
relationship \(\Q(C) = \Q(C^\star)\) from \Cref{lem:weighted-kleene-star-polyhedron}, 
leads to another of our main results.

\begin{theorem}\label{thm:polyhedral-cone-of-modification}
  Fix a DAG \(\dag\) with weights \(w\in\RR^E\).
  \begin{enumerate}
    \item The set of weights \(w'\in\RR^{E^\star}\) such that \(\Q(C) = \Q(C')\) is the pointed polyhedral cone \[
            \mathcal{C}(G,w) \defeq w^\star + \ray{e_{j\to i} \mid j\to i\in\dag^\star\setminus \dag^\flat}.
          \]
    \item The compactification of \(\mathcal{C}(G,w)\) inside \(\TT^{E^\star}\) with respect to its recession cone 
        is the set of all weights \(w'\) for any DAG \(\subdag=(V,E')\) satisfying 
    \(\dag^\flat_w\subseteq\subdag\subseteq\dag^\star\) such that \(\Q(C) = \Q(C')\).
    \item The set of \(w'\in\RR^E\) such that \(\Q(C) = \Q(C')\) is a face of the compactification \(\overline{\mathcal{C}(G,w)}\). 
  \end{enumerate}
  \begin{proof}\leavevmode
    \begin{enumerate}
        \item In the case that \(w'\) is chosen from the polyhedral cone \(\mathcal{C}(G,w)\)
        we can apply \Cref{lem:transitive-reduction-step} for every \(j\to i\notin\dag^\flat\)
        to conclude \(\Q(C') = \Q(C^\star)\) where \(C'\) is the weight matrix associated to \(w'\). 
        If on the other hand \(w'\in\RR^E\) such that \(\Q(C') = \Q(C)\), we know from 
        \Cref{lem:weighted-kleene-star-polyhedron} that additionally \(\Q(C') = \Q(C^\star) = \Q(C)\) must be satisfied. 
        If \(\{x_i - x_j = c'_{ij}\}\) is a facet-defining hyperplane of \(\Q(C)\), we must have \(c'_{ij} = c_{ij}^\star = c_{ij}\).
        On the other hand, if \(\{x_i - x_j = c'_{ij}\}\) does not define a facet, then \(c'_{ij} \geq c_{ij}^\star\).
        In other words, \(w'(j\to i) = w^\star(j\to i) + \varepsilon\) for some \(\varepsilon\in[0,\infty]\)
        which by \Cref{lem:transitive-reduction-step} is only possible for positive \(\varepsilon\) 
        if \(j\to i\notin G^\flat\).

        \item First, note that since \(\mathcal{C}(G,w)\) is a pointed polyhedral cone, it is equal to its own recession cone.
        In addition, the ray generators of \(\mathcal{C}(G,w)\) are given by coordinate directions \(e_{j\to i}\) for which 
        \(j\to i\notin G^\flat_w\). By \Cref{lem:compactifying-fan} this means that the compactification \(\overline{\mathcal{C}(G,w)}\) inside \(\TT^E\) 
        only intersects those strata indexed by subgraphs \(G^\flat_w\subseteq\subdag\subseteq G^\star\). 
        Now let \(w'\in\overline{\mathcal{C}(G,w)}\) be weights for such a subgraph.
        We know that \[
            w'(j\to i) = \begin{cases}
                w^\star(j\to i),& j\to i\in\dag^\flat,\\
                w^\star(j\to i) + \lambda_{j\to i}e_{j\to i},& j\to i\in H\setminus\dag^\flat_w.
            \end{cases}
        \] Let \(w''\in\RR^{E^\star}\) be weights on \(\dag^\star\) such that \(w''\vert_H = w'\) and 
        \(w''(j\to i) = w^\star(j\to i)\) otherwise. We may apply \Cref{cor:transitive-reduction-step} to obtain that 
        \(\Q(C') = \Q(C'') = \Q(C)\).

        If \(\Q(C') = \Q(C)\) for some graph \(\dag^\flat_w\subseteq\subdag\subseteq\dag^\star\), 
        then also \(\Q(C') = \Q(C^\star)\). Applying \Cref{lem:transitive-reduction-step} when \(w'(j\to i)>w^\star(j\to i)\)
        and \Cref{cor:transitive-reduction-step} when \(w'(j\to i) = w^\star(j\to i)\) gives \(w'\in\overline{\mathcal{C}(G,w)}\).

        \item This follows from the second point since \(G^\flat_w\subseteq G\subseteq G^\star\). \qedhere
    \end{enumerate}
  \end{proof}
\end{theorem}

\begin{theorem}\label{thm:transitive-reduction-minimal-facets}
    The weighted transitive reduction \(\dag^\flat_w\) is the unique minimal weighted DAG and \(w^\flat\) the unique weights
    such that \(\Q(C)\) and \(\Q(C^\flat)\) coincide.
    \begin{proof}
        First, we get \(\Q(C) = \Q(C^\flat)\) by applying \Cref{cor:transitive-reduction-step}
        for every edge \(j\to i\notin\dag^\flat_w\). \Cref{thm:polyhedral-cone-of-modification}.(2) implies that 
        \(\dag^\flat_w\) is the minimal graph for which this is possible.
        Assume that \(w'\) is another vector of weights for \(\dag^\flat_w\) and \(C'\) the corresponding weight matrix.
        If \(\Q(C') = \Q(C^\flat)\), we know by \Cref{lem:transitive-reduction-step} that for every \(j\to i\in\dag^\flat_w\)
        \(c'_{ij} = c^\flat_{ij}\). But since \(w'\) and \(w^\flat\) are weights for \(\dag^\flat_w\), 
        we already have \(w'=w^\flat\).
    \end{proof}
\end{theorem}

\begin{remark}\label{rem:stats-remark}
    With a view towards max-linear Bayesian networks, \Cref{thm:transitive-reduction-minimal-facets} describes exactly which weights of a network are identifiable given 
    sufficiently good data. Since \(\Q(C) = \Q(C^\star)\), sampling from a polytrope
    only allows us to recover \(C^\star\). 
    
    If \(C^\star = C^\flat\), this means that all original weights of \(C\) are identifiable.
    On the other hand, if there are edges \(j\to i\notin G^\flat\), it is not possible
    to estimate the corresponding weight \(c_{ij}\).
\end{remark}

As a consequence, the weighted transitive reduction \(C^\flat\) gives a finer correspondence between 
polytropes \(\Q(C)\) and the corresponding DAGs \(\dag\) than using the Kleene star, which would only
associate a transitive DAG to each \(C\). This also refines \Cref{thm:central-triangulations} 
in the following way.
\begin{corollary}\label{thm:central-triangulations-flat}
   If the regular subdivision \(\mathcal{S}\) of \(\fundamentalpolytope\)
   induced by \((0,w)\in\RR\times\RR^E\) is a central triangulation, then \(G^\flat_w = G\).
  \begin{proof}
    If \(C\in\TT^{n\times n}\) induces a central triangulation \(\subdivision\) on \(\fundamentalpolytope\), 
    for any \(j\to i\in G\), the edge between \(e_i - e_j\) and the origin occurs in \(\subdivision\).
    Any such edge is dual to a hyperplane \(\{x_i-x_j = c_{ij}\}\) and describes a facet of \(\Q(C)\). 
    But no hyperplane being redundant implies that \(C = C^\flat\), and therefore \(G^\flat_w = G\).
  \end{proof}
\end{corollary}

\begin{example}\label{ex:central-non-triangulation}
    \begin{figure}[tbp]







\begin{tikzpicture}[x  = {(0.9cm,-0.076cm)},
                    y  = {(-0.06cm,0.95cm)},
                    z  = {(-0.44cm,-0.29cm)},
                    scale = 2,
                    color = {lightgray}]

  \coordinate (v0_unnamed__1) at (0, 1, 0);
  \coordinate (v1_unnamed__1) at (0, 0, 1);
  \coordinate (v2_unnamed__1) at (-1, 1, 0);
  \coordinate (v3_unnamed__1) at (-1, 0, 1);
  \coordinate (v4_unnamed__1) at (0, 0, 0);

  \definecolor{vertexcolor_unnamed__1}{RGB}{ 186, 3, 252 }

  \tikzstyle{vertexstyle_unnamed__1_0} = [circle, scale=0.25, fill=vertexcolor_unnamed__1,]
  \tikzstyle{vertexstyle_unnamed__1_1} = [circle, scale=0.25, fill=vertexcolor_unnamed__1,]
  \tikzstyle{vertexstyle_unnamed__1_2} = [circle, scale=0.25, fill=vertexcolor_unnamed__1,label={[text=black, above left, align=left]:$e_3-e_2$},]
  \tikzstyle{vertexstyle_unnamed__1_3} = [circle, scale=0.25, fill=vertexcolor_unnamed__1,]
  \tikzstyle{vertexstyle_unnamed__1_4} = [circle, scale=0.25, fill=vertexcolor_unnamed__1,]

  \definecolor{facetcolor_unnamed__1}{RGB}{ 252, 186, 3 }

  \definecolor{edgecolor_unnamed__1}{rgb}{ 0 0 0 }
  \tikzstyle{facetstyle_unnamed__1} = [fill=facetcolor_unnamed__1, fill opacity=0.5, draw=edgecolor_unnamed__1, line width=1 pt, line cap=round, line join=round]

  \draw[facetstyle_unnamed__1] (v4_unnamed__1) -- (v3_unnamed__1) -- (v2_unnamed__1) -- (v4_unnamed__1) -- cycle;
  \draw[facetstyle_unnamed__1] (v1_unnamed__1) -- (v3_unnamed__1) -- (v4_unnamed__1) -- (v1_unnamed__1) -- cycle;
  \draw[facetstyle_unnamed__1] (v2_unnamed__1) -- (v0_unnamed__1) -- (v4_unnamed__1) -- (v2_unnamed__1) -- cycle;
  \draw[facetstyle_unnamed__1] (v2_unnamed__1) -- (v3_unnamed__1) -- (v1_unnamed__1) -- (v0_unnamed__1) -- (v2_unnamed__1) -- cycle;
  \draw[facetstyle_unnamed__1] (v0_unnamed__1) -- (v1_unnamed__1) -- (v4_unnamed__1) -- (v0_unnamed__1) -- cycle;

  \foreach \i in {3,1,2,4,0} {
    \node at (v\i_unnamed__1) [vertexstyle_unnamed__1_\i] {};
  }

  \coordinate (v0_unnamed__2) at (0, 0, 1);
  \coordinate (v1_unnamed__2) at (-1, 0, 1);
  \coordinate (v2_unnamed__2) at (0, -1, 1);
  \coordinate (v3_unnamed__2) at (0, 0, 0);

  \definecolor{vertexcolor_unnamed__2}{RGB}{ 186, 3, 252 }

  \tikzstyle{vertexstyle_unnamed__2_0} = [circle, scale=0.25, fill=vertexcolor_unnamed__2,]
  \tikzstyle{vertexstyle_unnamed__2_1} = [circle, scale=0.25, fill=vertexcolor_unnamed__2,label={[text=black, above left, align=left]:$e_4-e_2$},]
  \tikzstyle{vertexstyle_unnamed__2_2} = [circle, scale=0.25, fill=vertexcolor_unnamed__2,]
  \tikzstyle{vertexstyle_unnamed__2_3} = [circle, scale=0.25, fill=vertexcolor_unnamed__2,]

  \definecolor{facetcolor_unnamed__2}{RGB}{ 0,0,255 }

  \definecolor{edgecolor_unnamed__2}{rgb}{ 0 0 0 }
  \tikzstyle{facetstyle_unnamed__2} = [fill=facetcolor_unnamed__2, fill opacity=0.1, draw=edgecolor_unnamed__2, line width=1 pt, line cap=round, line join=round]

  \draw[facetstyle_unnamed__2] (v2_unnamed__2) -- (v1_unnamed__2) -- (v3_unnamed__2) -- (v2_unnamed__2) -- cycle;
  \draw[facetstyle_unnamed__2] (v0_unnamed__2) -- (v1_unnamed__2) -- (v2_unnamed__2) -- (v0_unnamed__2) -- cycle;
  \draw[facetstyle_unnamed__2] (v3_unnamed__2) -- (v1_unnamed__2) -- (v0_unnamed__2) -- (v3_unnamed__2) -- cycle;
  \draw[facetstyle_unnamed__2] (v2_unnamed__2) -- (v3_unnamed__2) -- (v0_unnamed__2) -- (v2_unnamed__2) -- cycle;

  \foreach \i in {1,2,0,3} {
    \node at (v\i_unnamed__2) [vertexstyle_unnamed__2_\i] {};
  }

  \coordinate (v0_unnamed__3) at (1, 0, 0);
  \coordinate (v1_unnamed__3) at (0, 1, 0);
  \coordinate (v2_unnamed__3) at (0, 0, 1);
  \coordinate (v3_unnamed__3) at (0, 0, 0);

  \definecolor{vertexcolor_unnamed__3}{RGB}{ 186, 3, 252 }

  \tikzstyle{vertexstyle_unnamed__3_0} = [circle, scale=0.25, fill=vertexcolor_unnamed__3,]
  \tikzstyle{vertexstyle_unnamed__3_1} = [circle, scale=0.25, fill=vertexcolor_unnamed__3,label={[text=black, above right, align=left]:$e_3-e_1$},]
  \tikzstyle{vertexstyle_unnamed__3_2} = [circle, scale=0.25, fill=vertexcolor_unnamed__3,label={[text=black, below left, align=left]:$e_4-e_1$},]
  \tikzstyle{vertexstyle_unnamed__3_3} = [circle, scale=0.25, fill=vertexcolor_unnamed__3,]

  \definecolor{facetcolor_unnamed__3}{RGB}{ 0,0,255 }

  \definecolor{edgecolor_unnamed__3}{rgb}{ 0 0 0 }
  \tikzstyle{facetstyle_unnamed__3} = [fill=facetcolor_unnamed__3, fill opacity=0.1, draw=edgecolor_unnamed__3, line width=1 pt, line cap=round, line join=round]

  \draw[facetstyle_unnamed__3] (v3_unnamed__3) -- (v2_unnamed__3) -- (v1_unnamed__3) -- (v3_unnamed__3) -- cycle;
  \draw[facetstyle_unnamed__3] (v0_unnamed__3) -- (v2_unnamed__3) -- (v3_unnamed__3) -- (v0_unnamed__3) -- cycle;
  \draw[facetstyle_unnamed__3] (v3_unnamed__3) -- (v1_unnamed__3) -- (v0_unnamed__3) -- (v3_unnamed__3) -- cycle;

   \node at (v3_unnamed__3) [vertexstyle_unnamed__3_3] {};

  \draw[facetstyle_unnamed__3] (v1_unnamed__3) -- (v2_unnamed__3) -- (v0_unnamed__3) -- (v1_unnamed__3) -- cycle;

  \foreach \i in {2,1,0} {
    \node at (v\i_unnamed__3) [vertexstyle_unnamed__3_\i] {};
  }

  \coordinate (v0_unnamed__4) at (1, 0, 0);
  \coordinate (v1_unnamed__4) at (0, 0, 1);
  \coordinate (v2_unnamed__4) at (0, -1, 1);
  \coordinate (v3_unnamed__4) at (0, 0, 0);

  \definecolor{vertexcolor_unnamed__4}{RGB}{ 186, 3, 252 }

  \tikzstyle{vertexstyle_unnamed__4_0} = [circle, scale=0.25, fill=vertexcolor_unnamed__4,label={[text=black, above right, align=left]:$e_2-e_1$},]
  \tikzstyle{vertexstyle_unnamed__4_1} = [circle, scale=0.25, fill=vertexcolor_unnamed__4,]
  \tikzstyle{vertexstyle_unnamed__4_2} = [circle, scale=0.25, fill=vertexcolor_unnamed__4,label={[text=black, below right, align=left]:$e_4-e_3$},]
  \tikzstyle{vertexstyle_unnamed__4_3} = [circle, scale=0.25, fill=vertexcolor_unnamed__4,label={[text=black, above right, align=left]:$0$},]

  \definecolor{facetcolor_unnamed__4}{RGB}{ 0,0,255 }

  \definecolor{edgecolor_unnamed__4}{rgb}{ 0 0 0 }
  \tikzstyle{facetstyle_unnamed__4} = [fill=facetcolor_unnamed__4, fill opacity=0.1, draw=edgecolor_unnamed__4, line width=1 pt, line cap=round, line join=round]

  \draw[facetstyle_unnamed__4] (v3_unnamed__4) -- (v2_unnamed__4) -- (v1_unnamed__4) -- (v3_unnamed__4) -- cycle;
  \draw[facetstyle_unnamed__4] (v0_unnamed__4) -- (v2_unnamed__4) -- (v3_unnamed__4) -- (v0_unnamed__4) -- cycle;
  \draw[facetstyle_unnamed__4] (v1_unnamed__4) -- (v2_unnamed__4) -- (v0_unnamed__4) -- (v1_unnamed__4) -- cycle;
  \draw[facetstyle_unnamed__4] (v3_unnamed__4) -- (v1_unnamed__4) -- (v0_unnamed__4) -- (v3_unnamed__4) -- cycle;

  \foreach \i in {2,1,3,0} {
    \node at (v\i_unnamed__4) [vertexstyle_unnamed__4_\i] {};
  }

\end{tikzpicture}
        \caption{The central regular subdivison from \Cref{ex:central-non-triangulation} with height
        function \((0,w) = (0,1,2,2,3,3,6)\). The weights \(w\) satisfy \(\dag^\flat_w = \dag\) but
        the subdivision is not a triangulation due to the orange cell.}
        \label{fig:central-non-triangulation}
    \end{figure}
    The following example shows that the converse of \Cref{thm:central-triangulations-flat}
    fails for $n\geq 4$.
    Let \(G = \kappa_4\) be the complete DAG on \(4\) nodes and choose as weights 
    \(w = (1,2,2,3,3,6)\). That is, the associated polytrope \(\Q(C)\) is given by the matrix \[
        C = \begin{pmatrix}
            0 & \infty & \infty & \infty \\
            1 & 0 & \infty & \infty \\
            2 & 3 & 0 & \infty \\
            2 & 3 & 6 & 0
        \end{pmatrix}.
    \] A straightforward calculation shows that \(C\) is idempotent, so \(C^\star = C\) and in particular
    \(\dag^\flat_w = \dag\). The corresponding subdivision of \(\fundamentalpolytope\) is shown in 
    \Cref{fig:central-non-triangulation}. This subdivision is central, but it contains the maximal cell 
    with vertices \[
        \{\, 0, e_3-e_1, e_3 - e_2, e_4 - e_1, e_4 - e_2 \,\}
    \] which is colored orange in \Cref{fig:central-non-triangulation}, and clearly not a simplex. 
    This means that the subdivision induced by \((0,w)\) is central, but not a triangulation.
\end{example}

\begin{remark}
    In practice, for \(C\in\TT^{n\times n}\) supported on a DAG, one can calculate the weighted transitive
    reduction by calculating \(C^\star\) and deleting every edge \(j\to i\)
    for which \(c_{ij} \neq c_{ij}^\star\). This amounts to the usual \(\bigo(n^3)\) 
    running time for calculating the Kleene star.
    
    Using Algorithm 2 of \citeauthor{Joswig.Schröter:2022} from \cite{Joswig.Schröter:2022},
    one can calculate from a given DAG \(G\) the possible weighted transitive reductions 
    with the corresponding inequalities. This may be computed in 
    \(\bigo(\lvert E\rvert^3\cdot\sigma + n^2)\) time where \(\sigma\) is the maximum number of possible
    shortest-path trees at any node.
\end{remark}

\subsection{Moduli space and its stratifications}

We discuss in the following the constructions of a moduli space \(\polytropemoduli{n}\) for 
\((n-1)\)-dimensional polytropes associated to DAGs.
For our construction, we use the secondary fan \(\sfan(\fundamentalpolytope)\) in \(\RR\times\RR^E\) 
of the fundamental polytope \(\fundamentalpolytope\) restricted to height functions \((0,w)\) where \(w\in\RR^E\)
are the original weights on \(\dag\).
We have two choices for how to write \(\polytropemoduli{n}\) as a collection of polyhedral fans.
\begin{enumerate}
    \item Classifying weights up to the Kleene stars of matrices \(C\in\TT^{n\times n}\) gives a covering
        by polyhedral fans indexed by transitive DAGs \(\dag\).
    \item Using the weighted transitive reduction \(C^\flat\), which results in a finer characterization in terms 
        of minimal non-redundant facet descriptions and the associated DAGs. This gives an atlas of charts 
        for \(\polytropemoduli{n}\) as a polyhedral space
\end{enumerate}

\begin{definition}
For any DAG \(\dag\), define the \emph{open polytrope region} \(\Kleene\) as the set of weights 
\(w\in\RR^{E}\) 
such that \(\dag\) with weights \(w\) agrees with its weighted transitive reduction \(\dag_w^\flat\). 
\end{definition}
From the definition of the weighted transitive reduction it follows that
the \emph{open polytrope region} \(\Kleene\) is an open polyhedral cone:
    {\small\[
        \Kleene = \left\{\, 
            w\in\RR^{E} \;\middle|\; w(j\to i) < \sum_{i = 0}^\ell w(v_i\to v_{i+1})
            \text{ for every path } 
            (j = v_0, v_1,\dots,v_\ell = i) \text { in } \dag
        \,\right\}.
    \]}

Restricting \(\sfan(\fundamentalpolytope)\) to \(0\times\Kleene\) gives a polyhedral fan
\(\localpolytropemoduli\) which we call the \emph{open stratum} for \(\dag\) and which by definition 
parametrizes polytropes \(\Q(C)\) where \(C^\flat\) is supported precisely on \(\dag\).

\begin{remark}
    One can also see the open polytrope region as a region in the tropical hypersurface arrangement given by the
    entries of the Kleene star for a generic matrix supported on \(\dag\).
\end{remark}

\begin{theorem}\label{thm:polytrope-moduli}\label{thm:fine-stratification}
    The polyhedral subspace of \(\TT^{E(\kappa_n)}\) given by the union of polyhedral fans \[
        \polytropemoduli{n} \defeq \coprod_{[\dag]} \localpolytropemoduli[\dag] \subset \TT^{E(\kappa_n)}
    \] where \([\dag]\) is taken over isomorphism classes of directed acyclic graphs on \(n\) nodes
    is a moduli space \(\polytropemoduli{n}\) for \((n-1)\)-dimensional polytropes associated to labeled DAGs.
    This means, every cone of \(\polytropemoduli{n}\) forms an equivalence class of polytropes up to tropical equivalence.

    \begin{proof}
        If \(P\) is a polytrope associated to a DAG, then there is a DAG \(\dag\) with weights \(w\in\RR^E\) and 
        weight matrix \(C\in\TT^{n\times n}\) such that \(P = \Q(C)\). By \Cref{thm:transitive-reduction-minimal-facets}, 
        \(\dag^\flat_w\) is the minimal DAG with unique weights \(w^\flat\) with this property. 
        As such, \(P\) corresponds to a point in the open stratum \(\localpolytropemoduli[\dag^\flat]\).

        Set-theoretically each open stratum \(\localpolytropemoduli\) is equal to the open polytrope region \(\Kleene\),
        which is an open subset of \(\RR^E\). Thus, the \(\localpolytropemoduli\) already provide an atlas of charts for 
        \(\polytropemoduli{n}\). Each chart is disjoint meaning there is nothing left to check for the transition maps.

        If \(\Q(C)\) and \(\Q(C')\) are two tropically equivalent polytropes such that \(C\) and \(C'\) are
        supported on DAGs \(\dag\) resp.\ \(\dag'\), we know there is a pair of permutations 
        \((\tau,\sigma)\in\SymGr_n\times\SymGr_n\) inducing an isomorphism on the poset of covectors.
        We may assume \wLOG\ that \(\Q(C)=\Q(C')\) up to a permutation of the coordinates. 
        Then, \(\Q(C^\flat)=\Q((C')^\flat)\) up to a permutation of the coordinates and in particular \(C^\flat = (C')^\flat\) 
        up to permutation of rows and columns. This means that \(w = w'\) up to permutation, 
        and the tropical equivalence of those two polytropes induces an isomorphism between \(G\) and \(G'\) preserving 
        the weights.

        On the other hand, if \(G\) and \(G'\) are isomorphic DAGs and again assume \wLOG\ that \(w = w'\) 
        up to permutation, then the associated weight matrices \(C\) resp.\ \(C'\) only differ by a 
        permutation of the rows and columns. But this gives a pair of permutations \((\tau,\sigma)\) to
        certify tropical equivalence.
    \end{proof}
\end{theorem}

If \(\dag\) is transitive, we can compare the weighted transitive reduction
to the Kleene star. In this case, the defining inequalities of \(\Kleene\) 
are just the strict triangle inequalities on \(\dag\). These are weights \(w\in\RR^E\) for which the 
associated weight matrix \(C\) satisfies both \(C^\star = C\) and \(C^\flat = C\).

The set of weights such that the associated weight matrix satisfies \(C=C^\star\) but not necessarily \(C^\flat = C\)
corresponds to weights satisfying the non-strict triangle inequalities and thus this set is the closure of \(\Kleene\) 
in \(\RR^E\).
Restricting a secondary fan \(\sfan(\fundamentalpolytope)\) for a transitive DAG \(\dag\) to the 
closed polytrope region \(0\times \overline{\Kleene}\) results in a polyhedral fan \(\closedlocalmoduli\) 
which we call the \emph{closed stratum} for \(\dag\).

The closed strata give a covering of \(\polytropemoduli{n}\) using less charts than in \Cref{thm:polytrope-moduli} in the following way.
\begin{proposition}\label{prop:open-closed-subfan}
    Let \(\dag\) be a transitive DAG and \(G^\flat\subseteq H\subseteq G^\star\). Then, \(\localpolytropemoduli[\subdag]\)
    can be identified with a subfan of \(\closedlocalmoduli\) supported on the boundary of \(\Kleene\).

    \begin{proof}
        For \(w\in\localpolytropemoduli[\subdag]\), define the inclusion of \(\varphi\colon \localpolytropemoduli[\subdag]\hookrightarrow\closedlocalmoduli\)
        by \(w\mapsto w^\star\). Since \(G^\flat\subseteq H\), \(w^\star\) will be a point in \(\RR^{E^\star}\) 
        and by construction \(w^\star\) satisfies the non-strict triangle inequalites, so \(w^\star\in\closedlocalmoduli\).

        Since \(\subdag\subsetneq\dag\), \(w^\star\) will satisfy the triangle equality \(w^\star(j\to i) = w^\star(j\to k) + w^\star(k\to i)\),
        for some edges \(j\to k,k\to i\in\subdag\) and \(j\to i\in\dag\), meaning that \(w^\star\) lies on the boundary of 
        \(\Kleene\).

        We need to show that cones of \(\localpolytropemoduli[\subdag]\) get mapped to cones in \(\closedlocalmoduli\).
        Note that the mapping \(w\mapsto w^\star\) defines an extension of the regular subdivision 
        \(\subdivision(\fundamentalpolytope, (0,w))\) to a regular subdivision 
        \(\subdivision(\fundamentalpolytope[\dag^\star], (0,w^\star))\). 
        
        Any new vertex \(e_i - e_j\) from \(\fundamentalpolytope[\dag^\star]\)
        is inserted at a height \(w^\star(j\to i)\) which is equal to the sum of heights of other vertices 
        already present in \(\fundamentalpolytope\). This means that any new vertex \(e_i - e_j\) only extends 
        maximal cells of \(\subdivision(\fundamentalpolytope, (0,w))\) without subdividing those any further.
        This means there is a one-to-one mapping of maximal cells of \(\subdivision(\fundamentalpolytope, (0,w))\)
        and \(\subdivision(\fundamentalpolytope[\dag^\star], (0,w^\star))\). But this means that two 
        points in the image of \(\varphi\) give equivalent regular subdivisions on \(\fundamentalpolytope[\dag^\star]\)
        if and only if points in their preimage give equivalent regular subdivisons of \(\fundamentalpolytope\).
    \end{proof}
\end{proposition}

\begin{figure}[tp]
    \centering
    \begin{tikzpicture}[x=4.5cm,y=4cm]
\node[draw, circle, inner sep=-7pt] (K3) at (0,0) {
    \begin{tikzcd}[row sep=.5em]
            1 \arrow[rd] \arrow[dd] & \\
            & 2 \arrow[ld] \\
            3
    \end{tikzcd}};
\node[draw, circle, inner sep=-7pt] (line) at (0,-1) {
    \begin{tikzcd}[row sep=.5em]
            1 \arrow[rd] & \\
            & 2 \arrow[ld] \\
            3
    \end{tikzcd}};
\node[draw, circle, inner sep=-7pt] (cherry) at (-1,-1) {
    \begin{tikzcd}[row sep=.5em]
            1 \arrow[rd] \arrow[dd] & \\
            & 2 \\
            3
    \end{tikzcd}};
\node[draw, circle, inner sep=-7pt] (collider) at (1,-1) {
    \begin{tikzcd}[row sep=.5em]
            1 \arrow[dd] & \\
            & 2 \arrow[ld] \\
            3
    \end{tikzcd}};
\node[draw, circle, inner sep=-7pt] (12) at (0,-2) {
    \begin{tikzcd}[row sep=.5em]
            1 \arrow[rd]  & \\
            & 2 \\
            3
    \end{tikzcd}};

\node[draw, circle, inner sep=-7pt] (empty) at (0,-3) {
    \begin{tikzcd}[row sep=.5em]
            1 & \\
            & 2 \\
            3
    \end{tikzcd}};

\coordinate[above left=of K3] (class1-1);
\coordinate[above right=of K3] (class1-2);
\coordinate[below right=of line] (class1-3);
\coordinate[below left=of line] (class1-4);


\begin{scope}[rounded corners=30, very thick, blue]
    \draw (class1-1) -- (class1-2) -- (class1-3) -- (class1-4) -- cycle;
\end{scope}

\begin{scope}[very thick]
\draw[-] (K3) -- (line);
\draw[-] (K3) -- (collider);
\draw[-] (K3) -- (cherry);

\draw[-] (12) -- (line);
\draw[-] (12) -- (line);
\draw[-] (12) -- (collider);
\draw[-] (12) -- (collider);
\draw[-] (12) -- (cherry);
\draw[-] (12) -- (cherry);

\draw[-] (empty) -- (12);
\draw[-] (empty) -- (12);
\draw[-] (empty) -- (12);
\end{scope}
\end{tikzpicture}
    \caption{The lattice of DAGs on \(n=3\) nodes up to isomorphism. 
    The single non-trivial closed stratum is circled.}
    \label{fig:graph-poset-3}
\end{figure}
\begin{figure}[bp]
    \centering
    \resizebox{.7\linewidth}{!}{\newsavebox\completedag
\begin{lrbox}{\completedag}
\adjustbox{scale=.75}{\begin{tikzcd}[row sep=.5em,blue!85]
  1  \arrow[dd] \arrow[rd] & \\
  & 2 \arrow[ld] \\
  3
\end{tikzcd}}
\end{lrbox}

\newsavebox\cherry
\begin{lrbox}{\cherry}
\adjustbox{scale=.75}{\begin{tikzcd}[row sep=.5em,magenta]
  1  \arrow[dd] \arrow[rd] & \\
  & 2 \\
  3
\end{tikzcd}}
\end{lrbox}

\newsavebox\chain
\begin{lrbox}{\chain}
\adjustbox{scale=.75}{\begin{tikzcd}[row sep=.5em,orange!70]
  1 \arrow[rd] & \\
  & 2 \arrow[ld] \\
  3
\end{tikzcd}}
\end{lrbox}

\newsavebox\collider
\begin{lrbox}{\collider}
\adjustbox{scale=.75}{\begin{tikzcd}[row sep=.5em,green!70!black]
    1  \arrow[dd] & \\
    & 2 \arrow[ld] \\
    3
\end{tikzcd}}
\end{lrbox}

\newsavebox\singleedge
\begin{lrbox}{\singleedge}
\adjustbox{scale=.75}{\begin{tikzcd}[row sep=.5em,black]
    1  \arrow[dd] & \\
    & 2 \\
    3
\end{tikzcd}}
\end{lrbox}
\begin{tikzpicture}[
    scale=4,
    roundnode/.style={circle, draw=black!85, fill=black!85, thin, scale=0.5},
    pseudonode/.style={circle, draw=black!50, fill=black!50, thin, scale=0.5},
    line width=3pt
]
    \coordinate (e1) at (1,1);
    \coordinate (e2) at (-1,0);
    \coordinate (e3) at (0,-1);

    \coordinate (i1) at (-1,-1);
    \coordinate (i2) at (1,0);
    \coordinate (i3) at (0,1);

    \fill[orange!10] (e1.center) -- (i3.center) -- (e2.center) -- (i1.center) -- cycle;
    \fill[blue!10] (e3.center) -- (i1.center) -- node[anchor=north west] {\usebox\completedag} (e1.center) -- (i2.center) -- cycle;
    

    \draw[ultra thick,orange!70] (i1) -- node[anchor=south east] {\usebox\chain} (e1) -- (i3) -- (e2);
    \draw[ultra thick,orange!70] (i1) -- (e2);

    \draw[ultra thick,green!70!black] (e2) -- node[anchor=east] {\usebox\collider} (i1) -- (e3);

    \draw[ultra thick,magenta] (e1) -- node[anchor=west] {\usebox\cherry} (i2) -- (e3);
    
    \node[roundnode, label={[label distance=1pt]30:$(0,\infty,\infty)$}] at (e1) {};
    \node[anchor=south east] at (e2) {\usebox\singleedge};
    \node[roundnode, label={[label distance=1pt]180:$(\infty,0,\infty)$}] at (e2) {};
    \node[roundnode, label={[label distance=1pt]-30:$(\infty,\infty,0)$}] at (e3) {};

    \node[pseudonode, draw=black!50!green!40, fill=black!50!green!40,] at (i1) {};
    \node[pseudonode, draw=black!50!magenta!40, fill=black!50!magenta!40,] at (i2) {};
    \node[pseudonode, draw=black!50!orange!40, fill=black!50!orange!40,] at (i3) {};
\end{tikzpicture}}
    \caption{The tropical projective plane \(\TTPP^2\) as the space of weights \(w\) 
    subdivided by combinatorial types.
    For each colored region, the corresponding weighted transitive reduction \(\dag^\flat_w\)
    is shown. For the vertices of \(\TTPP^2\), only one tropical vertex (in black) is marked since the 
    combinatorial type is the same for all tropical vertices. The interior of \(\TTPP^2\) corresponds to 
    weights on the complete DAG \(\kappa_3\), while the blue region is the open polytrope region \(\Kleene[\kappa_3]\) and the blue region together with the orange line segment is the 
    closure of \(\Kleene[\kappa_3]\).}
    \label{fig:param-space-3}
\end{figure}

\begin{corollary}\label{thm:coarse-stratification}
    The union of closed strata \(\closedlocalmoduli\) where \(\dag\) ranges over isomorphism classes of 
    transitive DAGs is a covering of the moduli space \(\polytropemoduli{n}\).

    \begin{proof}
        Note that any DAG is (non-uniquely) a subgraph of a transitive DAG. \Cref{prop:open-closed-subfan} then gives the result.
    \end{proof}
\end{corollary}

\begin{figure}[btp]
    \centering
    \begin{tikzpicture}[x=4.5cm,y=4cm]
\node[draw, circle, inner sep=-7pt] (K4) at (0,0) {
    \begin{tikzcd}[row sep=1.5em]
            1 \arrow[d] \arrow[rd] \arrow[r] & 2 \arrow[ld] \arrow[d]\\
            4 & 3 \arrow[l]
    \end{tikzcd}};

\node[draw, circle, inner sep=-7pt] (K4-13) at (-1,-1) {
    \begin{tikzcd}[row sep=1.5em]
            1 \arrow[d] \arrow[r] & 2 \arrow[ld] \arrow[d]\\
            4 & 3 \arrow[l]
    \end{tikzcd}};
\node[draw, circle, inner sep=-7pt] (K4-14) at (-0,-1) {
    \begin{tikzcd}[row sep=1.5em]
            1 \arrow[rd] \arrow[r] & 2 \arrow[ld] \arrow[d]\\
            4 & 3 \arrow[l]
    \end{tikzcd}};
\node[draw, circle, inner sep=-7pt] (K4-24) at (1,-1) {
    \begin{tikzcd}[row sep=1.5em]
            1 \arrow[d] \arrow[rd] \arrow[r] & 2 \arrow[d]\\
            4 & 3 \arrow[l]
    \end{tikzcd}};

\node[draw, circle, inner sep=-7pt] (K4-13-14) at (-1,-2) {
    \begin{tikzcd}[row sep=1.5em]
            1 \arrow[r] & 2 \arrow[ld] \arrow[d]\\
            4 & 3 \arrow[l]
    \end{tikzcd}};
\node[draw, circle, inner sep=-7pt] (K4-13-24) at (0,-2) {
    \begin{tikzcd}[row sep=1.5em]
            1 \arrow[d] \arrow[r] & 2 \arrow[d]\\
            4 & 3 \arrow[l]
    \end{tikzcd}};
\node[draw, circle, inner sep=-7pt] (K4-14-24) at (1,-2) {
    \begin{tikzcd}[row sep=1.5em]
            1 \arrow[rd] \arrow[r] & 2 \arrow[d]\\
            4 & 3 \arrow[l]
    \end{tikzcd}};
    
\node[draw, circle, inner sep=-7pt] (line) at (0,-3) {
    \begin{tikzcd}[row sep=1.5em]
            1 \arrow[r] & 2 \arrow[d]\\
            4 & 3 \arrow[l]
    \end{tikzcd}};

\begin{scope}[xshift=9cm,yshift=-4cm]
\node[draw, circle, inner sep=-7pt] (diamond) at (0,0) {
    \begin{tikzcd}[row sep=1.5em]
            1 \arrow[d] \arrow[rd] \arrow[r] & 2 \arrow[ld]\\
            4 & 3 \arrow[l]
    \end{tikzcd}};
 \node[draw, circle, inner sep=-7pt] (diamond-14) at (1,-1) {
    \begin{tikzcd}[row sep=1.5em]
            1 \arrow[rd] \arrow[r] & 2 \arrow[ld]\\
            4 & 3 \arrow[l]
    \end{tikzcd}};
\end{scope}

\begin{scope}[xshift=9cm,yshift=-8cm]
\node[draw, circle, inner sep=-7pt] (K4-12) at (0,0) {
    \begin{tikzcd}[row sep=1.5em]
            1 \arrow[d] \arrow[rd] & 2 \arrow[ld] \arrow[d]\\
            4 & 3 \arrow[l]
    \end{tikzcd}};
\node[draw, circle, inner sep=-7pt] (K4-12-14) at (0,-1) {
    \begin{tikzcd}[row sep=1.5em]
            1 \arrow[rd] & 2 \arrow[ld] \arrow[d]\\
            4 & 3 \arrow[l]
    \end{tikzcd}};
\node[draw, circle, inner sep=-7pt] (K4-12-23) at (1,-1) {
    \begin{tikzcd}[row sep=1.5em]
            1 \arrow[d] \arrow[rd] & 2 \arrow[ld]\\
            4 & 3 \arrow[l]
    \end{tikzcd}};
\node[draw, circle, inner sep=-7pt] (K4-12-14-23) at (1,-2) {
    \begin{tikzcd}[row sep=1.5em]
            1 \arrow[rd] & 2 \arrow[ld]\\
            4 & 3 \arrow[l]
    \end{tikzcd}};
\end{scope}

\begin{scope}[xshift=-9cm,yshift=-4cm]
\node[draw, circle, inner sep=-7pt] (K4-34) at (0,0) {
    \begin{tikzcd}[row sep=1.5em]
            1 \arrow[d] \arrow[rd] \arrow[r] & 2 \arrow[ld] \arrow[d]\\
            4 & 3 
    \end{tikzcd}};
\node[draw, circle, inner sep=-7pt] (K4-13-34) at (-1,-1) {
    \begin{tikzcd}[row sep=1.5em]
            1 \arrow[d] \arrow[r] & 2 \arrow[ld] \arrow[d]\\
            4 & 3 
    \end{tikzcd}};
\node[draw, circle, inner sep=-7pt] (K4-14-34) at (0,-1) {
    \begin{tikzcd}[row sep=1.5em]
            1 \arrow[rd] \arrow[r] & 2 \arrow[ld] \arrow[d]\\
            4 & 3 
    \end{tikzcd}};
\node[draw, circle, inner sep=-7pt] (K4-13-14-34) at (-1,-2) {
    \begin{tikzcd}[row sep=1.5em]
            1 \arrow[r] & 2 \arrow[ld] \arrow[d]\\
            4 & 3 
    \end{tikzcd}};
\end{scope}

\begin{scope}[xshift=4.5cm,yshift=-12cm]
\node[draw, circle, inner sep=-7pt] (K4-24-34) at (0,0) {
    \begin{tikzcd}[row sep=1.5em]
            1 \arrow[d] \arrow[rd] \arrow[r] & 2 \arrow[d]\\
            4 & 3 
    \end{tikzcd}};
\node[draw, circle, inner sep=-7pt] (K4-13-24-34) at (0,-1) {
    \begin{tikzcd}[row sep=1.5em]
            1 \arrow[d] \arrow[r] & 2 \arrow[d]\\
            4 & 3 
    \end{tikzcd}};
\end{scope}

\begin{scope}[xshift=-9cm,yshift=-12cm]
\node[draw, circle, inner sep=-7pt] (K4-12-34) at (0,0) {
    \begin{tikzcd}[row sep=1.5em]
            1 \arrow[d] \arrow[rd] & 2 \arrow[ld] \arrow[d]\\
            4 & 3 
    \end{tikzcd}};
\end{scope}


\begin{scope}[xshift=-4.5cm,yshift=-12cm]
\node[draw, circle, inner sep=-7pt] (K4-13-23) at (0,0) {
    \begin{tikzcd}[row sep=1.5em]
            1 \arrow[d] \arrow[r] \arrow[dr] & 2 \\
            4 & 3 \arrow[l]
    \end{tikzcd}};
\node[draw, circle, inner sep=-7pt] (K4-13-14-23) at (0,-1) {
    \begin{tikzcd}[row sep=1.5em]
            1 \arrow[rd] \arrow[r] & 2 \\
            4 & 3 \arrow[l]
    \end{tikzcd}};
\end{scope}

\begin{scope}[xshift=-9cm,yshift=-16cm]
\node[draw, circle, inner sep=-7pt] (14-24-34) at (0,0) {
    \begin{tikzcd}[row sep=1.5em]
            1 \arrow[d] \arrow[rd] \arrow[r] & 2\\
            4 & 3 
    \end{tikzcd}};
\end{scope}

\begin{scope}[xshift=-13.5cm,yshift=-16cm]
\node[draw, circle, inner sep=-7pt] (14-23-24) at (0,0) {
    \begin{tikzcd}[row sep=1.5em]
            1 \arrow[d] & 2 \arrow[d]\arrow[ld]\\
            4 & 3 
    \end{tikzcd}};
\end{scope}
\begin{scope}[xshift=9cm,yshift=-16cm]
\node[draw, circle, inner sep=-7pt] (14-24-34) at (0,0) {
    \begin{tikzcd}[row sep=1.5em]
            1 \arrow[d] & 2 \arrow[ld]\\
            4 & 3 \arrow[l]
    \end{tikzcd}};
\end{scope}

\begin{scope}[xshift=-9cm,yshift=-20cm]
\node[draw, circle, inner sep=-7pt] (14-23) at (0,0) {
    \begin{tikzcd}[row sep=1.5em]
            1 \arrow[d] & 2 \arrow[d]\\
            4 & 3 
    \end{tikzcd}};
\end{scope}

\begin{scope}[xshift=-4.5cm,yshift=-20cm]
\node[draw, circle, inner sep=-7pt] (12-13) at (0,0) {
    \begin{tikzcd}[row sep=1.5em]
            1 \arrow[r]\arrow[rd] & 2 \\
            4 & 3 
    \end{tikzcd}};
\end{scope}

\begin{scope}[xshift=0cm,yshift=-24cm]
\node[draw, circle, inner sep=-7pt] (12) at (0,0) {
    \begin{tikzcd}[row sep=1.5em]
            1 \arrow[r] & 2 \\
            4 & 3 
    \end{tikzcd}};
\end{scope}

\begin{scope}[xshift=0cm,yshift=-28cm]
\node[draw, circle, inner sep=-7pt] (empty) at (0,0) {
    \begin{tikzcd}[row sep=1.5em]
            1 & 2 \\
            4 & 3 
    \end{tikzcd}};
\end{scope}

\begin{scope}[xshift=4.5cm,yshift=-20cm]
\node[draw, circle, inner sep=-7pt] (13-23) at (0,0) {
    \begin{tikzcd}[row sep=1.5em]
            1 \arrow[rd] & 2 \arrow[d]\\
            4 & 3
    \end{tikzcd}};
\end{scope}

\begin{scope}[xshift=0cm,yshift=-16cm]
\node[draw, circle, inner sep=-7pt] (12-13-23) at (0,0) {
    \begin{tikzcd}[row sep=1.5em]
            1 \arrow[r] \arrow[rd] & 2 \arrow[d]\\
            4 & 3
    \end{tikzcd}};
\node[draw, circle, inner sep=-7pt] (12-23) at (0,-1) {
    \begin{tikzcd}[row sep=1.5em]
            1 \arrow[r] & 2 \arrow[d]\\
            4 & 3
    \end{tikzcd}};
\end{scope}

\begin{scope}[very thick]
\draw[-] (K4) -- (K4-13);
\draw[-] (K4) -- (K4-14);
\draw[-] (K4) -- (K4-24);

\draw[-] (K4-13) -- (K4-13-14);
\draw[-] (K4-13) -- (K4-13-24);

\draw[-] (K4-14) -- (K4-14-24);
\draw[-] (K4-14) -- (K4-13-14);

\draw[-] (K4-24) -- (K4-14-24);
\draw[-] (K4-24) -- (K4-13-24);

\draw[-] (K4-13-14) -- (line);
\draw[-] (K4-13-24) -- (line);
\draw[-] (K4-14-24) -- (line);

\draw[-] (diamond) -- (diamond-14);

\draw[-] (K4-12) -- (K4-12-14);
\draw[-] (K4-12) -- (K4-12-23);
\draw[-] (K4-12-14) -- (K4-12-14-23);
\draw[-] (K4-12-23) -- (K4-12-14-23);

\draw[-] (K4-34) -- (K4-14-34);
\draw[-] (K4-34) -- (K4-13-34);
\draw[-] (K4-14-34) -- (K4-13-14-34);
\draw[-] (K4-13-34) -- (K4-13-14-34);


\draw[-] (K4-13-23) -- (K4-13-14-23);
\draw[-] (K4-24-34) -- (K4-13-24-34);

\draw[-] (12-13-23) -- (12-23);
\end{scope}
\end{tikzpicture}
    \caption{Representation of the lattice of DAGs on \(n=4\) nodes up isomorphism.
    Only the covering relations between DAGs associated with each closed stratum are shown.
    Equivalently, each connected component represents a closed local moduli space 
    \(\closedlocalmoduli[\dag]\) according to \Cref{thm:coarse-stratification}.}
    \label{fig:graph-poset-4}
\end{figure}
\begin{example}
    The graphs indexing the open strata are shown in \Cref{fig:graph-poset-3}.
    There is only one non-trivial closed stratum which is circled. 
    This is precisely the stratum that we discussed in \Cref{ex:transitive-reduction}.

    For comparison, \Cref{fig:param-space-3} shows the sketch of a subdivision of points \(w\in\TTPP^2\)
    according to the combinatorial type of polytrope \(P = \Q(C)\). In the case \(n=3\),
    combinatorial types are in one-to-one correspondence with types of weighted transitive reductions 
    \(\dag^\flat_w\). The blue circled interval of graphs from \Cref{fig:graph-poset-3} corresponds
    to the interior of \(\TTPP^2\) in \Cref{fig:param-space-3}.

    The secondary fan of \(\fundamentalpolytope[\kappa_3]\) for the directed triangle 
    has two maximal cones given by the two halfspaces of the hyperplane \[
        c = a + b.
    \] The maximal cones correspond to the inequalities from \Cref{ex:transitive-reduction}.
    The restriction of \(\sfan(\fundamentalpolytope[\kappa_3])\) to the open polytrope region
    is precisely the halfspace given by \(c<a+b\), which corresponds to the blue region in \Cref{fig:param-space-3}.
    
    Furthermore, the restriction to the closed polytrope region includes the cone
    where \(c = a+b\), which corresponds to \(\localpolytropemoduli[\subdag]\)
    for the path graph \(H\).
    This is the orange line in the interior of \(\TTPP^2\) of \Cref{fig:param-space-3}.

    The orange sector extends from the interior of \(\TTPP^2\) to a boundary component with \(c = \infty\). 
    This is due to \Cref{thm:polyhedral-cone-of-modification} which states that for every point in the orange region,
    there is a ray extending in the \(c\)-coordinate direction along which \(\Q(C)\) is constant.
    Then, the compactification of this ray adds a point on the orange boundary component.
    
    Any other graph \(\dag\) on three nodes does not correspond to the closed stratum \(\closedlocalmoduli[\kappa_3]\), 
    so there is no copy of \(\localpolytropemoduli\) in \(\closedlocalmoduli[\kappa_3]\).
\end{example}

For \(n=4\), only a partial representation is shown in \Cref{fig:graph-poset-4} for the sake of readability. 
Only the covering relations for the graphs associated to each closed stratum are shown, making each stratum visible as a 
connected component. Also, each connected component shows the structure of each closed stratum
\(\closedlocalmoduli[\dag^\star]\) in terms of the open strata \(\localpolytropemoduli[\subdag]\).

\section{Enumerating combinatorial types}\label{sec:enumeration}

Enumerating the tropical combinatorial types of polytropes coming from DAGs on \(n\) nodes
amounts to checking secondary fans ranging over a suitable set of DAGs. Due to the two different stratifications,
we have the choice between two different procedures to enumerate.

In light of \Cref{thm:coarse-stratification}, we can check the subdivisions of \(\rootpolytope\) for each 
transitive DAG \(\dag\) on \(n\) nodes. This gives the possible tropical combinatorial types for tropical polytopes.
Then, using \Cref{thm:polytrope-cell}, we need to take the face figure of the subdivision contributing to the polytrope 
and eliminate duplicates. It can be seen in \Cref{tab:enumeration} that this reduces the number of
DAGs to be checked.

One disadvantage of this approach is that we need to enumerate all regular subdivisions and 
can only filter duplicates afterwards. \Cref{tab:enumeration} gives a lower bound for the amount of subdivisions that need to be
stored.
It has been noted by \citeauthorandcite{Joswig.Schröter:2019} that the root polytope \(\rootpolytope\) 
is at most \(2n\)-dimensional with \(\bigo(n^2)\) vertices and lattice volume bounded by \(\binom{2n-2}{n-1}\).
Passing to the fundamental polytope \(\fundamentalpolytope\) reduces the dimensionality and the number of vertices.

Using \Cref{thm:fine-stratification} we can instead enumerate the regular central subdivisions of 
\(\fundamentalpolytope\) for every graph \(\dag\) without filtering, which alleviates the storage issue. 
When further restricting to regular central \emph{triangulations}, we can take advantage of the software \verb|mptopcom| 
by \citeauthorandcite{JJK:2018} that has been optimized for parallel enumeration of central triangulations.

Passing from subdivisons to triangulations restricts our focus to \emph{generic} combinatorial types in the sense
that the corresponding polytrope has the maximal number of faces for each dimension among all polytropes associated 
to a fixed graph \(\dag\). In other words, the generic combinatorial types correspond to the maximal cones of each
\(\localpolytropemoduli\). In \cite{Tran:2016}, \citeauthor{Tran:2016} also enumerated the maximal cones of a specific 
Gröbner fan to count the number combinatorial types that are generic for the complete graph \(K_n\) for \(n\leq 5\),
which she termed as maximal polytropes.

To enumerate central regular triangulations, we generated the input for \verb|mptopcom| using Julia~1.10.9 and {OSCAR}~1.4.1 
and used \verb|mptopcom| 1.4 for the enumeration.
This way, we can take advantage of the techniques implemented in \texttt{mptopcom} to avoid checking unnecessary triangulations 
and to enumerate orbits of triangulations.

The result of the enumerations are shown in \Cref{tab:enumeration}, with the number of DAGs that had to be enumerated
listed for each method. The code for performing the enumeration can be found under 
\href{https://zenodo.org/records/14166530}{https://zenodo.org/records/14166530}.

Comparing with the results from \citeauthorandcite{Tran:2016}, there are \(1013\) combinatorial types of bounded polytropes 
for \(n=4\) and \(27248\) combinatorial types of generic bounded polytropes for \(n=5\). 
In both cases we count significantly fewer generic combinatorial types of polytropes associated to DAGs.

\begin{table}[bt]
   \caption{Number of generic combinatorial types for polytropes associated to max-linear Bayesian networks in 
   \(\mathbb{R}^n/\mathbb{R}\mathbf{1}\) and the number of checked DAGs}
   \begin{tabular}{lccccccc}
       \toprule
       $n$ & 1 & 2 & 3 &  4 &    5 &     6 &  7 \\\midrule
        \#triangulations & 
             1 & 2 & 6 & 32 &  512 & 65240 &  276733462 \\
        \#DAGs &
             1 & 2 & 6 & 31 &  302 &  5984 &  243668  \\
        \#transitive DAGs &
             1 & 2 & 5 & 16 &   63 &   318 &  2045  \\
       \bottomrule
   \end{tabular}
   \label{tab:enumeration}
\end{table}

\section{Outlook}
In this work we have given different notions of equivalence for polytropes coming from max-linear Bayesian networks.
For a polytrope \(\Q(C)\) its minimal non-redundant facet description corresponds to a weight matrix \(C^\flat\) 
which is associated to the weighted transitive reduction. If \(\dag\) is the underlying DAG of \(C\),
the weighted transitive reduction has a natural interpretation in terms of 
those edges which are essential for realizing the shortest paths in \(\dag\). 
This provides an equivalence relation by labeling a polytrope with a DAG which is more specific than 
the underlying graph of the Kleene star. This graph is always a transitive DAG  
and does not see the redundancy of hyperplanes for defining facets.

We focused on DAGs because they serve as a description for max-linear Bayesian networks, and our construction of the weighted transitive reduction makes use of the absence of cycles.
With respect to previous work on polytropes which focused on the complete graph, a characterization 
for all other graphs which contain cycles is missing.

\begin{question}
    How does one characterize a minimal non-redundant facet description for arbitrary digraphs, 
    thus generalizing the notion of the weighted transitive reduction to graphs with cycles?
\end{question}

Seeing a polytrope \(\Q(C)\) as a subcomplex of the covector decomposition \(\CovDec(C)\) allowed us to 
characterize tropical equivalence in terms of regular subdivisions, meaning the secondary fan 
is the object of choice for constructing a parameter space. While this idea has been known in tropical geometry 
for some time, the construction with the secondary fan requires a fixed support graph \(\supportgraph[C]\). 
Hence, the main challenge when classifying all polytropes of max-linear Bayesian networks consists in 
fitting together the information from all necessary secondary fans for the possible support graphs.

Related work by \citeauthorandcite{OPS:2019} studied in depth \emph{local Dressians} for a given matroid.
In particular, they investigated the structure given by the secondary fan of the matroid polytope and
by the tropical variety cut out by three-term Plücker relations, showing these structures agree.
From this, the next step for max-linear Bayesian networks is to investigate the structure coming from
the Plücker relations.

Another point of view is that since matrices associated to networks are lower-triangular, we can consider \(\Q(C)\)
as a flag of tropical polytopes. Tropical flag varieties have been studied by
\citeauthorandcite{BEZ:2021} who found incidence-Plücker relations in addition to
the usual Plücker relations.
\begin{question}
    In what way do our local moduli spaces \(\localpolytropemoduli\), 
    the local Dressians for the cycle matroid of \(G\) 
    and the tropical varieties of complete flags relate to each other?
\end{question}

Finally, in terms of statistical applications, knowing that the support of a max-linear Bayesian network is a polytrope gives us additional information useful for modeling extreme data.
As stated in \Cref{rem:stats-remark}, only the weights associated to edges in the weighted transitive reduction 
can theoretically be recovered from samples. We also would like to know what further statistical consequences can be obtained from understanding the combinatorial structure, ideally shedding light on questions considered by \cite{amendola2021markov,amendola2022conditional}.

\begin{question}
   To what extent does the combinatorics of the polytrope associated to a max-linear Bayesian network reveal information about the conditional independence statements that hold among the variables (and viceversa)? 
\end{question}

\singlespace
\printbibliography
\end{document}